\documentclass[a4paper]{amsart}
\usepackage{graphicx} % Required for inserting images

\title{Growing spines: ad infinitum et ad infinitesimalia}

\author[B.~Boissonneau]{Blaise Boissonneau}
\address{Blaise Boissonneau \\ Heinrich Heine University Düsseldorf, Faculty of Mathematics and Natural Sciences, Universitätsstr.~1, 40225 Düsseldorf, Germany. \url{blaise.boissonneau@hhu.de}}

\author[A.~De Mase]{Anna De Mase}
\address{Anna De Mase \\ Università degli Studi della Campania "Luigi Vanvitelli", Dipartimento di Matematica e Fisica, viale A. Lincoln 5, 81100 Caserta, Italy and Università degli Studi Roma Tre, Dipartimento di Matematica e Fisica, via della vasca navale 84, 00146 Rome, Italy and Zukunftskolleg, University of Konstanz, Universitätsstraße 10,  78464 Konstanz, Germany. \url{anna.de-mase@uni-konstanz.de}}
\thanks{ADM was partially supported by GNSAGA (INdAM, Italian National Institute of High Mathematics) - Project ``Model theory and applications to algebra and combinatorics", and PRIN 2022 ``Models, sets and classification".}

\author[F.~Jahnke]{Franziska Jahnke}
\address{Franziska Jahnke \\ Institute for Logic, Language and Computation, University of Amsterdam, 
Science Park 107,
1098 XG Amsterdam, The Netherlands 
and Institute for Mathematical Logic and Foundations, Department of Mathematics and Computer Science,
University of M\"unster,
Einsteinstraße 62,
48149 M\"unster, Germany. \url{franziska.jahnke@uni-muenster.de}}
\thanks{FJ was supported by the Deutsche Forschungsgemeinschaft (DFG, German Research Foundation) under Germany's Excellence Strategy EXC 2044-390685587, Mathematics M\"unster: Dynamics-Geometry-Structure, as well as by a Fellowship from the Daimler and Benz Foundation.}

\author[P.~Touchard]{Pierre Touchard}

\address{Pierre Touchard \\ Institut f\"ur Algebra, Technische Universit\"at Dresden, 01062 Dresden, Germany and KU Leuven, Department of Mathematics, B-3001 Leuven, Belgium. \url{pierre.touchard@tudresden.de}}

\thanks{PT was partially supported by KU Leuven IF C16/23/010.}
\thanks{All authors were partially funded by the Deutsche Forschungsgemeinschaft (DFG, German Research Foundation) under Germany Excellence Strategy-EXC-2047/1 - 390685813 via the Hausdorff Research Institute for Mathematics.}
\date{\today}
\usepackage[T1]{fontenc} 
\usepackage[utf8]{inputenc}
\usepackage{marvosym}
\usepackage{amsmath}
\usepackage{scalerel}
\usepackage{amsthm}
\usepackage{amsfonts}
\usepackage{amssymb}
\usepackage{comment}
\usepackage{graphicx}
\usepackage{slashed} 
\usepackage{braket}
\usepackage{hyperref}
\usepackage[nameinlink]{cleveref}
\usepackage{tikz-cd}
\usepackage{adjustbox}
\usepackage{euscript}
\usepackage{thmtools}
\usepackage{thm-restate}
\usepackage{xparse}
\usepackage{mathtools}
\usetikzlibrary{decorations.pathreplacing,angles,quotes}

\usepackage[inline]{enumitem}

\hypersetup{
  colorlinks   = true, %Colours links instead of ugly boxes
  urlcolor     = {blue!50!black}, %Colour for external hyperlinks
  linkcolor    = {blue!50!black}, %Colour of internal links
  citecolor   = {red!50!black} %Colour of citations
}
\def\restriction#1#2{\mathchoice
              {\setbox1\hbox{${\displaystyle #1}_{\scriptstyle #2}$}
              \restrictionaux{#1}{#2}}
              {\setbox1\hbox{${\textstyle #1}_{\scriptstyle #2}$}
              \restrictionaux{#1}{#2}}
              {\setbox1\hbox{${\scriptstyle #1}_{\scriptscriptstyle #2}$}
              \restrictionaux{#1}{#2}}
              {\setbox1\hbox{${\scriptscriptstyle #1}_{\scriptscriptstyle #2}$}
              \restrictionaux{#1}{#2}}}
\def\restrictionaux#1#2{{#1\,\smash{\vrule height .8\ht1 depth .85\dp1}}_{\,#2}}

\DeclareMathOperator{\val}{val}

\DeclareMathOperator{\res}{res}

\DeclareMathOperator{\Sp}{Sp}

\DeclareMathOperator{\Th}{Th}

\DeclareMathOperator{\CVX}{Conv}
\DeclareMathOperator{\Char}{char}

\makeatletter
\DeclareRobustCommand\bigop[1]{%
  \mathop{\vphantom{\sum}\mathpalette\bigop@{#1}}\slimits@
}
\newcommand{\bigop@}[2]{%
  \vcenter{%
    \sbox\z@{$#1\sum$}%
    \hbox{\resizebox{\ifx#1\displaystyle.9\fi\dimexpr\ht\z@+\dp\z@}{!}{$\m@th#2$}}%
  }%
}
\makeatother

\newcommand{\hprod}{\DOTSB\bigop{\mathrm{H}}}

\newcommand{\Ccal}{\ensuremath{\mathcal{C}}}

\newcommand{\Lcal}{\ensuremath{\mathcal{L}}}

\newcommand{\Ucal}{\ensuremath{\mathcal{U}}}

\newcommand{\Nbb}{\ensuremath{\mathbb{N}}}

\newcommand{\Qbb}{\ensuremath{\mathbb{Q}}}
\newcommand{\Rbb}{\ensuremath{\mathbb{R}}}

\newcommand{\Zbb}{\ensuremath{\mathbb{Z}}}

\newcommand{\Ztwo}{\ensuremath{\mathbb{Z}_{(2)}}}
\newcommand{\Zthree}{\ensuremath{\mathbb{Z}_{(3)}}}

\newcommand\cla\EuScript
\let\lang\CMcal

\newcommand{\lL}{\lang{L}}

\usepackage{ulem}
\usepackage{todonotes}
\setuptodonotes{noinline}

\colorlet{PierreColor}{cyan}

\colorlet{BlaiseColor}{violet!80!white}

\colorlet{AnnaColor}{green!80!yellow}

\colorlet{FranziColor}{yellow}

\NewDocumentCommand{\sset}{mg}{\left\{#1\IfNoValueF{#2}{\mid#2}\right\}}

\theoremstyle{plain}
\newtheorem{theorem}{Theorem}[section]
\newtheorem{thmx}{Theorem}
\newtheorem{lemma}[theorem]{Lemma}
\newtheorem{proposition}[theorem]{Proposition}
\newtheorem{fact}[theorem]{Fact}
\newtheorem{corollary}[theorem]{Corollary}

\theoremstyle{definition}
\newtheorem{definition}[theorem]{Definition}

\newtheorem*{notation*}{Notation}
\newtheorem{question}[theorem]{Question}
\newtheorem*{question*}{Question}
\newtheorem{example}[theorem]{Example}

\theoremstyle{remark}

\newtheorem{remark}[theorem]{Remark}
\newtheorem*{remark*}{Remark}

\newtheorem{observation}[theorem]{Observation}
\newtheorem{claim}{Claim}
\newtheorem*{claim*}{Claim}
\newtheorem{poc}{Proof of Claim}

\definecolor{black}{rgb}{0,0,0}

\colorlet{savedColor}{.}

\crefname{subsection}{subsection}{subsections}
\crefname{subsubsection}{subsubsection}{subsubsections}

\newcommand{\Lval}{\ensuremath{\mathcal L_{\mathrm{val}}}}
\newcommand{\Lring}{\ensuremath{\mathcal L_{\mathrm{ring}}}}
\newcommand{\Loag}{\ensuremath{\mathcal L_{\mathrm{OAG}}}}

\let\vphi\varphi

\let\substruct\preccurlyeq
\let\superstruct\succcurlyeq

\newcommand{\emtpyset}{\emptyset}
\newcommand{\Gdiv}{\ensuremath{G^{\mathrm{div}}}}

\renewcommand{\phi}{\varphi}

\begin{document}

\begin{abstract}

We prove that for every ordered abelian group $G$ there exists a non-trivial ordered abelian group $H$ such that $G\substruct H\oplus G$ with the lexicographic order, and give a first-order characterization of ordered abelian group $G$ such that $G\substruct G\oplus H$ for some non-trivial $H$. 
We apply this to characterize which ordered abelian groups (respectively fields)
ensure that any henselian valuation with said value group (respectively 
residue field)
is definable in the language of rings. This answers a question of  Krapp, 
Kuhlmann, and Link.
\end{abstract}
\maketitle

%auto-ignore
\section{Introduction}
This paper investigates \emph{automatic definability} of henselian valuations, a notion motivated by questions raised by Krapp, Kuhlmann, and Link \cite{KKL}.
%, and aims to provide a complete overview.%:which ordered abelian groups $\Gamma$ have the property that if $(K,v)$ is a henselian valued field with value group $\Gamma$, then $v$ is definable (respectively $\emptyset$-definable) in the language of rings? Similarly, which fields $k$ have the property that if $(K,v)$ is a henselian valued field with residue field $k$, then $v$ is definable (respectively $\emptyset$-definable) in the language of rings?
%We provide answers to these questions. 
Along the way, we study a range of structures (namely valued fields, ordered abelian groups and coloured linear orders) and provide characterizations 
of ``augmentability'' in these settings. 
%This connection, studied in detail in \Cref{sec:def-and-aug}, establishes a framework for studying definability problems, which give rise to interesting definability results. %not only allows us to address the motivating questions but overall gives rise to many definability results.

\subsection{Automatic Definability}

We call a valuation $v$ on a field $K$ \emph{definable} if the valuation ring $\mathcal O_v$ is a definable set in the language of rings.

It is a well-known fact that if the value group of a henselian valuation is $\mathbb Z$, then the valuation is definable, and even $\emptyset$-definable, see \cite{Robinson} and \cite{Ax71}. In \cite{Hong}, it is shown more generally that the valuation is definable if the value group is $p$-regular and not $p$-divisible. This inspires the following terminology: %We say that these valuations are \emph{automatically definable} by properties of the value group, and by analogy, we define:

\begin{definition}\ 
\begin{itemize}
    \item An ordered abelian group $G$ has \emph{automatic definability} (in the class of henselian valued fields) if in any henselian valued field $(K,v)$ with value group $vK=G$, the valuation is definable. We say that the valuation $v$ is definable \emph{by properties of the value group}.
    \item A field $k$ has \emph{automatic definability} if in any henselian valued field $(K,v)$ with residue field $Kv=k$, the valuation is definable. We say that the valuation $v$ is definable \emph{by properties of the residue field}.
\end{itemize}
\end{definition}

 We define similarly \emph{automatic $\emptyset$-definability} by requiring valuations to be $\emptyset$-definable. By Beth's Definability Theorem, 
 one sees that automatic $\emptyset$-definability in fact coincides with 
 uniform $\emptyset$-definability of the valuation across the set of Henselian valued fields with a fixed residue or a fixed value group. In this terminology, 
 the questions raised by 
 Krapp, Kuhlmann, and Link are: 
\begin{question*}[{\cite[Question 5.1]{KKL}}]
Which ordered abelian groups and which fields have automatic ($\emptyset$-)definability?
\end{question*}

We give a full characterization for ordered abelian groups, and for fields of characteristic $0$. More precisely, we show:

\begin{thmx}\label{thmA}
    Let $k$ be a field of characteristic $0$. The following are equivalent:
    \begin{enumerate}
        \item $k$ is t-henselian, i.e it is elementary equivalent to a field $k'$ that admits a non-trivial henselian valuation,
        \item for some valued field $(K,v)$ with residue field $k$, the valuation $v$ is not definable (not even over parameters), i.e., $k$ does not have automatic definability, 
        \item for some valued field $(K,v)$ with residue $k$, the valuation $v$ is not $\emptyset$-definable, i.e., $k$ does not have automatic $\emptyset$-definability.
    \end{enumerate}
\end{thmx}

\noindent This is proven in \Cref{def_fieldwise}. As we discuss in \Cref{sec:AutomaticDefinability}, this theorem does not hold in positive characteristic.
Similar characterizations for existential (respectively universal) definability of valuations were shown by Anscombe and Fehm \cite{AF17}.

%Note that the theorem does not hold for imperfect residue fields $k$: If $k$ is imperfect and admits no Galois extensions of degree divisible by $p$ (e.g., $k=\mathbb{F}_p(t)^\mathrm{sep}$), then any henselian valuation with residue field $k$ is definable by the Scanlon-trick (Proposition 3.6 in Jahnke: When does NIP transfer...?). Note that 
%$\mathbb{F}_p(t)^\mathrm{sep}$ is t-henselian (and in fact even henselian).\blaise{maybe remove this and leave it as a remark whereever it is stated in section 3}

For ordered abelian groups, the characterization is formulated in terms of \emph{augmentability}: an ordered abelian group $G$ is called 
\emph{weakly} (resp.~\emph{strongly}) \emph{augmentable by infinitesimals} if there is a non-trivial ordered abelian group $H$ such that $G \equiv G \oplus H$ (resp.~$G \substruct G \oplus H$). We show:

 \begin{thmx}\label{thmB}
    Let $G$ be an ordered abelian group. The following are equivalent:
    \begin{enumerate}
        \item $G$ is weakly augmentable by infinitesimals,
        \item for some valued field of $(K,v)$ with value group $G$, the valuation $v$ is not $\emptyset$-definable, i.e., $G$ does not have automatic $\emptyset$-definability.
    \end{enumerate}
    Moreover, the following are equivalent:
    \begin{enumerate}
        \item $G$ is strongly augmentable by infinitesimals,
        \item $G$ has no minimum positive element and is closed in its divisible hull $\Gdiv$,
        \item for some valued field of $(K,v)$ with value group $G$, the valuation $v$ is not definable (not even over parameters), i.e., $G$ does not have automatic definability,
    \end{enumerate}
 \end{thmx}
 
This is the combination of \Cref{weakABC,str-abc,Gdiv}. This theorem in turn 
gives rise
to a new question:
\begin{question} \label{que:2} Which ordered abelian groups are augmentable by infinitesimals, both weakly or strongly? 
\end{question}
Note that while \Cref{Gdiv} provides a characterization of ordered abelian groups having automatic definability
without referring to strong augmentability, its proof
relies on the classification of strongly augmentable ordered abelian groups.

%One deduces from this the aforementioned result of 
%\cite{Hong}, that is, 
%As a consequence, we recover Hong's result \cite{Hong} that
%regular non-divisible ordered abelian groups have automatic $\emptyset$-definability, since they are not weakly augmentable.

\subsection{Augmentable ordered abelian groups}
Ordered abelian groups are much-studied objects in model theory. By the results of Schmitt in his Habilitationsschrift \cite{Sch82}, and Cluckers-Halupczok \cite{CH11}, many first-order properties of ordered abelian groups can be reduced to corresponding properties of their \emph{spines}, which are chains of uniformly definable convex subgroups. %However, due to the highly technical nature of the framework, most recent works restrict to subclasses, e.g., ordered abelian groups with finite spines, and results applying to the class of all ordered abelian groups are rare.  
To answer \Cref{que:2}, we will be using the formalism of spines of Schmitt, and building on the impressive work of Delon and Lucas in \cite{DelonLucas}. We prove:

\begin{thmx}\label{thmD}
An ordered abelian group $G$ is strongly augmentable by infinitesimals if and only if for any $n\geqslant 2$, there exists no $g\in G$ with $F_n(g)=\{0\}$.

Additionally, $G$ is weakly augmentable by infinitesimals if and only if there exists $k\geqslant 2$ such that for any $n\geqslant 2$ with $k|n$, the $n$-Spine of $G$ is weakly augmentable on the right, that is, there's a non-empty coloured linear order $\Sigma_n$ such that $\Sp_n(G)\equiv\Sp_n(G)+\Sigma_n$. 
\end{thmx}

\noindent \Cref{thmD} is proven as \Cref{thm:StronglyAugmentableInfinitesimalsOAGs} and \Cref{thm:WeakAugmentableInfinitesimalsOags}. For an introduction to the formalism of spines, see
\Cref{subsec:lang-of-spines}. We also characterize augmentability by infinites, by proving that it holds for every ordered abelian groups:

\begin{thmx}
 Let $G$ be an ordered abelian group. Then $G$ is strongly augmentable by infinites, that is, there always exists a non-trivial ordered abelian group $H$ such that $G\substruct H\oplus G$.
\end{thmx}

\noindent This is \Cref{thm:strongleftaugmentabilityOAG}. It is a key part of our proof of \Cref{thmA}. One may inquire when a group is ``augmentable'' by a specific group $H$. We do not explore this question in depth, but we characterize ordered abelian groups $G$ with $G\substruct G\oplus\mathbb Q$ in \Cref{DOAG-aug} and those with $G\substruct\mathbb Q\oplus G$ in \Cref{cor_div}.

Once again, the characterization of weak augmentability of $G$ is phrased
in terms of weak augmentability of its spine. Thus, the following question
remains:
\begin{question}\label{que:3}
    When is a coloured linear order augmentable?
\end{question}

\subsection{Augmentable linear orders}
Model theory of linear orders is also a topic on his own, one of the main reference being \cite{Rub74}. To the best of our knowledge \Cref{que:3} was not explicitly studied before, although it is close in spirit to other
results (see, e.g. \cite{EP}). We show the following:

\begin{thmx}
 Let $X$ be a coloured linear order. Call $X$ strongly augmentable (on the right) if there exists $Y$ non empty such that $X\substruct X+Y$, and weakly augmentable if $X\equiv X+Y$. Then:
 \begin{enumerate}
 \item $X$ is strongly augmentable if and only if it does not have a maximal element,
 \item $X$ is weakly augmentable if and only if it is elementary equivalent to something of the form $B+\sum_{\omega^*}L$, where $B$ and $L$ are non-empty coloured linear orders.
 \end{enumerate}
\end{thmx}

This is  \Cref{prop:ChainAugmentable} and \Cref{lem:aug-mo} below. The first characterization (existence of a maximal element) is finitely axiomatisable, while the second one is axiomatisable but not finitely axiomatisable, as we explain in \Cref{prop:theory-of-weak-aug-lo}.
With these results, we conclude our study.

\bigskip 

To summarize,
%We present our results in as straightforward a manner as possible. 
we group our results by theme rather than ordering them by logical implications, therefore, some proofs rely on results presented later in the paper. 

First, we present in \Cref{sec:def-and-aug} a wide variety of results tying together augmentability and definability. This can serve as a motivation for the rest of the paper, and we also believe some results to be of interest in and of themselves.

Then, we characterize augmentability of ordered abelian groups in \Cref{sec:aug-oag}. This relies on the machinery of spines, which we recall in detail in \Cref{subsec:lang-of-spines} and \Cref{subsec:theory-of-spines}. With the help of these results, we characterize automatic definability for ordered abelian groups and for fields in \Cref{sec:AutomaticDefinability}. Finally, we study augmentable linear orders in \Cref{sec:aug-lo}, furthering our understanding of augmentability.

\subsection{Historical remarks on the model theory of ordered abelian groups}
%\franzi{I have not changed anything here yet, but I plan to}
Ever since the fundamental results of Robinson and Zakon \cite{RoZa60}, the
model theory of ordered abelian groups has been studied intensely.
Schmitt wrote his Habilitationsschrift \cite{Sch82} in 1982, which collects most of his work on ordered abelian groups, including his celebrated collaboration with Gurevich, published two years later (\cite{GS84}). It builds on and completes Gurevich’s work \cite{Gur64,Gur77}. This Habilitationsschrift 
is a milestone in the model theory of ordered abelian groups: the principle of reducing the theory of ordered abelian groups to theories of coloured chains is clearly exposed, the theory of spines is fully axiomatised, etc. In turn, general properties of coloured chains can then be lifted to ordered abelian groups. Among other applications, let us mention a cancellation theorem proved simultaneously by Delon-Lucas \cite{DelonLucas} and Giraudet \cite{Gir88}:  if $A$ and $B$ are ordered abelian groups and $A^n \equiv B^n$, then already $A \equiv B$. 
In fact, in \cite{DelonLucas}, Delon and Lucas enriched the theory of Schmitt and proved various striking results on ordered abelian groups. We will recall some of them in this paper.  Surprisingly, their paper appears to go largely unnoticed (we found two citations, only one recorded on MathSciNet). 

Consecutively, Delon and Farré published in 1996 a work on almost real closed fields \cite{DF96}, where again the theory of Schmitt is used and enriched. Remarkably, a double reduction, similar to the one we use in the present paper (from Henselian valued fields to ordered abelian groups, then from ordered abelian groups to coloured chains), structures their paper. To our knowledge, this is the first time that Schmitt's work is used in order to study definability in henselian valued fields.   

Outside of these authors (Delon, Farré, Lucas), few seem to cite Schmitt and Gurevich's work in a non-superficial way. Even in \cite{CH11}, where the authors reformulate Schmitt's reduction in more modern terms (with notions such as, e.g., relative quantifier elimination), Schmitt's Habilitationsschrift is only cited for a syntactical comparison. 

Contemporary to Schmitt, we find other major works in the subject, with authors such as Weispfenning \cite{Wei81, Wei86} and Giraudet \cite{Gir88}, mentioned earlier. Both use different sets of techniques.

Modern works are numerous, and we only cite a few: distality in ordered abelian groups is partially studied in \cite{ACGZ20}, a complete classification of strongly dependent ordered abelian groups can be found in \cite{HH19A}, and elimination of imaginaries is partially treated in \cite{Vic22}. Due to the nature of these questions, or by a deliberate choice to restrict them, authors tend only to study ordered abelian groups with finite spines, eluding completely the difficulties that Schmitt's theory is meant to address.
It is perhaps not an exaggeration to say that Schmitt’s work has been largely forgotten or overlooked in modern model theory. This is certainly due to several factors, such as the nature of the Habilitationsschrift (non-published, difficult to find) and the (unavoidable) technical nature of the statements. This contrasts with the striking efficiency of this theory, to which \cite{DelonLucas} is perhaps the best illustration.

\subsection{Notation}

\subsubsection{Valued fields}
We usually denote a valued field $(K,v)$, where $v$ is the valuation map. We typically denote by $\mathcal{O}_v$ its valuation ring, by $vK$ its value group and by $Kv$ and its residue field. The first-order language of valued fields is $\Lval\coloneq\Lring\cup\{\mathcal O\}$, where $\Lring\coloneqq \{+,\cdot~,0,1 \}$ is the language of rings, and $\mathcal O$ is a unary predicate interpreted
as the valuation ring.

If $k$ is a field and $\Gamma$ an ordered abelian group, $k((\Gamma))$ denotes the Hahn series field, that is, the field of formal power series $\sum_{g\in \Gamma} a_g t^g$, where $a_g\in K$, with well-ordered support. This field is equipped with the $t$-adic valuation ``$\val$'' sending $\sum_{g\in \Gamma} a_g t^g$ to $\min\sset{g\in \Gamma}{a_g\neq0}$. Note that for $G$ and $H$ ordered abelian groups, $k((G\oplus H))\cong k((H))((G))$, where $\oplus$ is the lexicographic sum -- see below.

\subsubsection{Ordered abelian groups}
The first-order language of ordered abelian groups is $\Loag=\{0,+,<\}$. We denote by $\oplus$ the lexicographic sum of ordered abelian groups, that is, given $I$ an ordered index set and for each $i\in I$ an ordered abelian group $G_i$, $\bigoplus_{i\in I} G_i$ is the set of sequences $(g_i)_{i\in I}$, with each $g_i\in G_i$, of finite support (i.e. $g_i\neq 0$ for finitely many $i$), equipped with coordinate-wise addition and with the lexicographic order: $(g_i)_{i\in I}<(h_i)_{i\in I}$ if and only if there is $j\in I$ such that $g_j<h_j$ and $g_i=h_i$ for all $i<j$. Note that $H$ (or formally, $\{0\}\times H$) is a convex subgroup of $G\oplus H$.

We denote by $\hprod$ the Hahn product of ordered abelian groups, that is, $\hprod_{i\in I}G_i$ is the set of sequences with well-ordered support, coordinate-wise product and lexicographic order. If $g$ is an element of an ordered abelian group, we write $\vert g \vert $ for the absolute value of $g$, that is, $\vert g \vert =\max\{g,-g\}$. If $S$ is a subset of $G$, we let $\langle S \rangle$ denotes the group generated by $S$, and $\CVX_G\langle S \rangle$ or simply $\CVX\langle S \rangle$ the convex hull in $G$ of the group generated by $S$.

\subsubsection{Linear orders}

We often work with coloured linear orders, that is, in a language with a binary relation ``$<$'' and an arbitrary number of unary predicates. In a fixed such language, given an ordered index set $I$ and for each $i\in I$ a linear order $X_i$, we denote by $\sum_{i\in I} X_i$ the disjoint union of all the $X_i$, equipped with the ordering $x<y$ if and only if either $x\in X_i$, $y\in X_j$ and $i<j$, or $x,y\in X_i$ and $X_i\vDash x<y$. If we work with colours, for any unary predicate $P$ in the language and $x\in \sum_{i\in I} X_i$, we colour $x$ by $P$ if and only if $x\in X_i$ and $X_i\vDash P(x)$. 
If all $X_i$ are identical copies of the same linear order $X$, then the lexicographic sum $\sum_{i\in I} X_i$ is often denoted $Y\cdot X$ or $\sum_Y X$.

We denote simply by $n$ the linear order with $n$ points (usually all monochromatic). We write $\omega$ for the order-type of the natural numbers $\mathbb N$, $\omega^*$ for the set $\omega$ equipped with the reversed order, so $\omega^*\vDash x<y$ if and only if $\omega\vDash x>y$, $\zeta$ for the order type of the integers $\mathbb Z$, and $\eta$ for the order type of the rationals $\mathbb Q$.

\subsection*{Acknowledgments}
Our gratitude goes to Salma Kuhlmann for her comments on a previous version of this work and for pointing us to the work of Delon and Lucas. 
In our first ArXiv preprint, we were only able to show that all ordered abelian groups are augmentable by infinites -- essentially reproving one of the main theorems of \cite{DelonLucas}.
Building further on the
ideas of Delon and Lucas, we were able to obtain characterizations for 
(weak) augmentability by infinitesimals, completing the picture.

In no particular order, we thank Artem Chernikov, Philip Dittmann, Arno Fehm and Immanuel Halupczok for their remarks during presentations of this project. We thank 
Sylvy Anscombe for her generous comments and technical support. 
We deeply thank David Gonzalez and Doug Cenzer for all the interesting discussions on linear orders, which nourished \Cref{sec:aug-lo}. 
Finally, we give some special thanks to Françoise Delon for various discussions, and to Peter Schmitt for sending us a copy of his Habilitationsschrift.

Part of this project was carried out at the Hausdorff Research Institute for Mathematics (HIM), to which we extend our gratitude.
%https://www.mathematics.uni-bonn.de/him/service/copy_of_him-faq#Going}

%All authors were partially funded by the Deutsche Forschungsgemeinschaft (DFG, German Research Foundation) under Germany ́s Excellence Strategy – EXC-2047/1 – 390685813}
%auto-ignore
\section{Definability and Augmentability}\label{sec:def-and-aug}

This section establishes many results showing links between definability of certain sets and the notion of ``augmentability'' in an ordered or valued structure. Most importantly, we study the definability of 
valuation rings in henselian valued fields of equicharacteristic 0.  This is, in some sense, folklore material, although we presented it with a new perspective. This will be extremely useful when we investigate notions of \emph{automatic definability} (\Cref{sec:AutomaticDefinability}).

%In this setting, we analyse elementary embeddings of the residue field and combine them with the augmentability by infinitesimals of the value group to give a necessary and sufficient condition for the valuation to be ($\emptyset$)-definable.

We mostly rely on the following classical result:
\begin{fact}[Beth's definability theorem, {\cite[Thm.~6.6.4]{Hodges}}]
    Let $\Lcal$, $\Lcal'$ be first-order languages with $\Lcal\subseteq\Lcal'$. Let $T$ be a theory in $\Lcal'$ and $\phi(\bar{x})$ a formula of $\Lcal'$. Then the following are equivalent:
    \begin{itemize}
        \item If $A$ and $B$ are models of $T$ such that $\restriction{A}{\Lcal}=\restriction{B}{\Lcal}$, then $A\models\phi(\bar{a})$ if and only if $B\models\phi(\bar{a})$, for all tuples $\bar{a}\in A$;
        \item $\phi(\bar{x})$ is equivalent modulo $T$ to a formula $\psi(\bar{x})$ of $\Lcal$.
    \end{itemize} \label{fact:beth}
\end{fact}

Beth's definability theorem is central to our arguments, and it will be used on various occasions throughout the paper. Its strength is to turn notions of definability into more malleable notions. As we shall see, in equicharacteristic 0 henselian valued fields, definability of the valuation can be phrased in terms of \emph{non-co-augmentability} (see  \Cref{def:WeakAugmentableVF,def:strongaugmentableVF}). This equivalence, in turn, enables us to apply the Ax–Kochen/Ershov machinery and explore various refinements and variants of the problem of definability.

We will also give two more applications of Beth's definability theorem, by studying definability of convex subgroups of ordered abelian groups and of end segments of coloured linear orders.

\subsection{Augmentable Henselian Valued fields}

Recall that a valuation $v$ on a field $K$ is \emph{definable} if there is
an $\mathcal{L}_\textrm{ring}$-formula $\psi(x)$, possibly using parameters from $K$, such that
$\psi(K)=\mathcal{O}_v$ holds. If $v$ is definable via a formula $\psi$ that requires 
no parameters from $K$, we say it is \emph{$\emptyset$-definable}.

A direct application of Beth's definability theorem for valued fields gives the following result: 
\begin{corollary} \label{cor:beth}
    Let $(K,v)$ be a valued field. Then, the following are equivalent:
    \begin{itemize}
        \item $v$ is \emph{not} $\emptyset$-definable in the language of ring,
        \item there is a field $K^*$ with two distinct valuations
        $u$ and $w$ such that $(K,v) \equiv (K^*,u) \equiv (K^*,w)$.\footnote{Such a valued field $(K,v)$ was called \emph{self}-similar in \cite[Definition 5.2]{AJ18}.}      
    \end{itemize}
    \end{corollary}
\begin{proof}
    The equivalence follows from Beth's definability theorem applied with the language $\lL=\Lring$, $\lL'=\Lval$, the theory  $T=\Th_{\lL'}(K,v)$ and where $\phi(x)$ is the $\Lval$-formula asserting $x\in\mathcal O$.
\end{proof}

\begin{example}
    Consider the $x$-adic valuation $v_x$ and the $y$-adic valuation $v_y$ on the field $K=\Rbb(x,y)$;  then by symmetry, $(K,v_x) \equiv (K,v_y)$ and neither valuation $v_x$ nor $v_y$ is $\emptyset$-definable.  
\end{example}

    Similarly, we can characterize definability over parameters: 
    
\begin{corollary}\label{cor:beth2}
     Let $(K,v)$ be a valued field.  
    Then the following are equivalent:
    \begin{enumerate}
        \item $v$ is \emph{not} definable (over parameters) in the language of rings,
        \item there is a field embedding $K\subseteq K^*$ in an extension $K^*$ and two distinct valuations $w,u$
        on $K^*$ such that $(K,v)\substruct(K^*,u)$ and $(K,v)\substruct(K^*,w)$. %and $(K,v) \substruct (K^*,w)$.
%        \item there is a field $K^*$ with two distinct \emph{comparable} valuations
%        $u$ and $w$ such that $(K,v) \substruct (K^*,u)$ and $(K,v) \substruct (K^*,w)$.
%Only needed for co-augmentability
    \end{enumerate}
\end{corollary}
\begin{proof} The equivalence %$(1) \Longleftrightarrow (2)$
follows once again from Beth's definability theorem, applied this time with the language $\lL=\Lring(K)$ (with constant symbols for each element in $K$), 
$\lL'=\Lval(K)$,  the theory $T=\Th_{\lL'}(K,v)$ and where $\phi(x)$ is once again the $\Lval$-formula asserting $x\in \mathcal{O}$.
\end{proof}

We now want to study the analogues of \Cref{cor:beth,cor:beth2} in the henselian setting, but where we require in addition the 
valuations $u$ and $w$ on $K^*$ to be comparable. In this case, one will
be a coarsening of the other (and hence corresponds to a convex subgroup). Thus, we recall some facts about convex subgroups of ordered
abelian groups and (non-)comparability of henselian valuations.
%Recall that a divisible ordered abelian group has no nontrivial proper definable
% subgroups. 
\begin{fact} [e.g. {\cite[Corollary 3.3.38]{AvdDvdH} }]\label{fact:PureSaturatedExactSequencesSplit}
    Consider an exact sequence of abelian groups:
    \[0 \rightarrow A \rightarrow B \rightarrow C \rightarrow 0.\]
    Assume that $A$ is pure in $B$ (i.e. all elements $a\in A$ which are $n$-divisible in $B$ are $n$-divisible in $A$). Assume also that $A$ is $\aleph_1$-saturated. 
    Then the exact sequence splits:
    $B\cong A \times C$.
\end{fact}

This gives the following well-known consequence:
\begin{corollary} \label{cor:OagSaturatedExactSequenceSplit}
    Let $(A, B)$ be a $\aleph_{1}$-saturated pair of ordered abelian groups where $A$ is a convex subgroup of $B$. Then $B\cong B/A\oplus A$.
    
    In particular, if $G$ is an arbitrary ordered abelian group and $C \leq G$
    a convex subgroup, we have $G\equiv G/C\oplus C$.
    \end{corollary}
    The first statement is immediate (see, e.g. \cite[Exercise 2.33]{marker18}). The second one is obtained by taking an $\aleph_1$-saturated extension after naming the convex subgroup with a predicate. 

Next, we recall the following facts from \cite{EP05} concerning the collection of all henselian valuations existing on a fixed field $K$, including the definition of the canonical henselian valuation. Notice that all fields admit at least one henselian valuation since the trivial valuation is henselian. Consider the following sets \begin{align*}H_{1}(K)&=\{v\,\vert\,v\,\text{is henselian},\,Kv\,\text{is not separably closed}\},\\H_{2}(K)&=\{v\,\vert\,v\,\text{is henselian},\,Kv\,\text{is separably closed}\}.\end{align*}
Two valuations $v,w$ are said to be \emph{comparable} if $\mathcal O_{v}\subseteq\mathcal O_{w}$ or $\mathcal O_{w}\subseteq\mathcal O_{v}$ holds for their 
valuation rings $\mathcal O_{v}$ resp.~$\mathcal O_{w}$. If $\mathcal O_{v}\subseteq\mathcal O_{w}$, then we say  that $w$ is \emph{coarser} than $v$ (or a \emph{coarsening} of $v$) and that $v$ is \emph{finer} than $w$ (or a \emph{refinement} of $w$).

\begin{fact}[{\cite[Theorem 4.4.2]{EP05}}]\label{fact:canval}
    Let $K$ be a field. Then the henselian valuations in $H_{1}(K)$ are linearly ordered by inclusion. If $H_{2}(K)\neq\emptyset$, then there is a valuation $v_{K}$ in $H_{2}(K)$ which is the coarsest valuation in $H_{2}(K)$, and it is finer than any valuation in $H_{1}(K)$. If $H_{2}(K)=\emptyset$, then there is a valuation $v_{K}$ in $H_{1}(K)$ which is the finest valuation in $H_{1}(K)$. \label{fact:tree}
\end{fact}
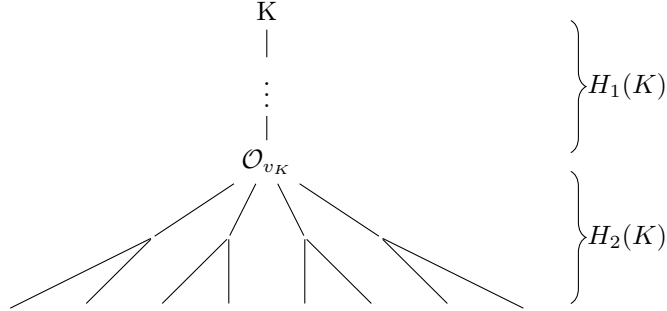
\begin{center}
    
\begin{tikzpicture}[align=center, sibling distance=1cm, level distance=1cm, node distance=2cm]
  \node(K) {K}
      child { node {\vdots}
	child { node {$\mathcal{O}_{v_K}$}
	  child{ node[inner sep=0]{}
	    child{ node{} }
	    child{ node{} }
	    child{ node{} edge from parent[draw=none]}
	    child{ node{} edge from parent[draw=none]}
	    child{ node{} edge from parent[draw=none]}
	  }
	  child{ node[inner sep=0]{}
	    child{ node{} edge from parent[draw=none]}
	    child{ node{} }
	    child{ node{} }
	    child{ node{} edge from parent[draw=none]}
	    child{ node{} edge from parent[draw=none]}
	  }
	  child{ node[inner sep=0]{}
	    child{ node{} edge from parent[draw=none]}
	    child{ node{} edge from parent[draw=none]}
	    child{ node{} }
	    child{ node{} }
	    child{ node{} edge from parent[draw=none]}
	  }
	  child{ node[inner sep=0]{}
	    child{ node{} edge from parent[draw=none]}
	    child{ node{} edge from parent[draw=none]}
	    child{ node{} edge from parent[draw=none]}
	    child{ node{} }
	    child{ node{} }
	  }
	}
      };
    \node(1) [right of=K] {};
    \node(top) [right of=1]   {};
    \node(mid) [below of=top] {};
    \node(bot) [below of=mid] {};
    
    \draw[decorate,decoration={brace,amplitude=0.2cm}] (top) to node[right] {\;$H_1(K)$} (mid);
    \draw[decorate,decoration={brace,amplitude=0.2cm}] (mid) to node[right] {\;$H_2(K)$} (bot);    
\end{tikzpicture}

\end{center}
In either case, the valuation $v_{K}$ is called the \emph{canonical valuation} of $K$, which is then comparable to any other henselian valuation on $K$. We will usually denote by $Kv_{K}$ and $v_{K}K$ the residue field and the value group of the canonical valuation.

Thus, the two valuations in the statement of \Cref{cor:beth} are always comparable if we assume one of the residue fields to be non-separably closed. If we restrict our study to equicharacteristic $0$ however, this last assumption is not necessary. Indeed, we can refine \Cref{cor:beth} as follows:

\begin{proposition}\label{prop:comparablevaluations0Definability}
    Let $(K,v)$ be a henselian valued field, with $\mathrm{char}(Kv)=0$. %non-separably closed. 
    Then the following are equivalent:
    \begin{enumerate}
        \item $v$ is \emph{not} $\emptyset$-definable in the language of rings,
        \item there is a field $K^*$ with two distinct \emph{comparable} valuations
        $u$ and $w$ such that $(K,v) \equiv (K^*,u) \equiv (K^*,w)$.
    \end{enumerate}
\end{proposition}
\begin{proof}
The direction (2) implies (1) follows immediately from \Cref{cor:beth}.
To see that (1) implies (2), note that \Cref{cor:beth} implies that if $v$ is not 
$\emptyset$-definable in $\mathcal{L}_\mathrm{ring}$, 
there is a field $K^*$ with two distinct valuations
        $u$ and $w$ such that $(K,v) \equiv (K^*,u) \equiv (K^*,w)$.
        In case $w$ and $u$ are already comparable, we are done.
%        If $Kv$ is not separably closed then the same holds for $K^*u$ and $K^*w$, hence $u$
%        and $w$ are comparable by \Cref{fact:tree}.
If they are not comparable, both their residue fields are algebraically closed
%If $Kv$ is not separably closed, then neither are $K^*w$ and $K^*u$, and so $w$ and $u$ are comparable 
by \Cref{fact:tree} --
%If $Kv$ is separably closed, it is necessarily also algebraically closed, and hence so are $K^*w$
%and $K^*u$.  Thus, we may assume that they are not comparable. In particular -- as in the proof of \Cref{prop:comparablevaluations0Definability} -- 
neither of them is equal to the canonical henselian valuation $v_{K^*}$ on $K^*$, and the residue field of $v_{K^*}$ is itself algebraically closed. 
In particular, $\Gamma^* =uK^*$ has a nontrivial divisible convex subgroup, since $u$ induces a nontrivial henselian valuation on the residue field of $v_{K^*}$ (see, e.g. \cite[p.45]{EP05}).
Hence, $\Gamma^*$ is elementarily equivalent to 
$\Gamma^* \oplus \mathbb{Q}$ (for a self-contained proof, see \Cref{DOAG-aug}).
%By \Cref{fact:EmbedsHahnSeries} and the elementary 
By the Ax-Kochen/Ershov principle (\cite[Theorem 6.17]{Hils-mtvf}), 
$(K, v)$ is elementarily equivalent to the field 
$Kv((\mathbb{Q}))(({\Gamma}^*))$ with two distinct but comparable valuations: one being the power series valuation with 
value group $\Gamma^* \oplus \mathbb{Q}$ and residue field $Kv$,
the other the power series valuation with value group $\Gamma^*$ and residue field $Kv((\mathbb{Q}))$ (note that the residue fields are elementary
equivalent as both are algebraically closed of the same characteristic).
\end{proof}

We show next that a similar result holds with parameters allowed. In the proof, we use the following fact (cf.~\cite{vdD14} for a proof).
\begin{fact}\label{fact:EmbedsHahnSeries}
    Let $(K,v)$ be a $\aleph_1$-saturated henselian valued field of equicharacteristic $0$. Then $(K,v)$ elementary embeds in the Hahn series $Kv((vK))$.
\end{fact}

\begin{proposition}\label{prop:comparablevaluationsDefinability}
    Let $(K,v)$ be a henselian valued field of equicharacteristic $0$.  
    Then the following are equivalent:
    \begin{enumerate}
        \item $v$ is \emph{not} definable (over parameters) in the language of rings,
        \item there is a field $K^*$ with two distinct \emph{comparable} valuations
        $u$ and $w$ such that $(K,v) \substruct (K^*,u)$ and $(K,v) \substruct (K^*,w)$.
    \end{enumerate}
\end{proposition}
\begin{proof}
    As before, the direction (2) implies (1) follows already from \Cref{cor:beth2}. It remains to show that (1) implies (2). The proof is 
    a variant of that of \Cref{prop:comparablevaluations0Definability}
Assume that $(K,v)$ is henselian with $\mathrm{char}(Kv)=0$ and that $v$ is not
$\mathcal{L}_\mathrm{ring}(K)$-definable.
 By \Cref{cor:beth2}, there is a field $K^*$ with two distinct valuations
$u$ and $w$ such that $(K,v) \substruct (K^*,u)$ and $(K,v) \substruct (K^*,w)$.
If $w$ and $u$ are already comparable, we are done.
If they are not comparable, both their residue fields are algebraically closed
%If $Kv$ is not separably closed, then neither are $K^*w$ and $K^*u$, and so $w$ and $u$ are comparable 
by \Cref{fact:tree} --
%If $Kv$ is separably closed, it is necessarily also algebraically closed, and hence so are $K^*w$
%and $K^*u$.  Thus, we may assume that they are not comparable. In particular -- as in the proof of \Cref{prop:comparablevaluations0Definability} -- 
neither of them is equal to the canonical henselian valuation $v_{K^*}$ on $K^*$. 
In particular, $\Gamma^* =uK^*$ has a nontrivial divisible convex subgroup
and hence elementarily embeds into $\Gamma^* \oplus \mathbb{Q}$ (for a self-contained proof, see \Cref{DOAG-aug}).
By \Cref{fact:EmbedsHahnSeries} and the elementary embedding version of the 
Ax-Kochen/Ershov principle (\cite[Theorem 6.17]{Hils-mtvf}), 
$(K, v)$ now embeds elementarily into the
Hahn product $K^*u((\mathbb{Q}))(({\Gamma}^*))$ with two different but comparable valuations: one embedding is with respect to the power series valuation with 
value group $\Gamma^* \oplus \mathbb{Q}$ and residue field $K^*u$,
the other with the valuation with value group $\Gamma^*$ and residue field $K^*u((\mathbb{Q}))$ (note that the embeddings of residue fields is elementary
as both are algebraically closed of the same characteristic).
%
%We consider now the elementary extension $ (K^{**},u^{**}):= k^*(({\Gamma}^*\oplus \Qbb))$. 
%We have a canonical elementary extension $(K^{*},u^{*}) \substruct (K^{**},u^{**})$. It is an elementary extension since $\Gamma^*$ has a nontrivial divisible convex subgroup, therefore ${\Gamma}^* \substruct {\Gamma}^*\oplus \Qbb$.
%On $K^{**}$, we can also consider a proper proper coarsening $v^{**}$ of $u^{**}$ with value group ${\Gamma}^*$ and residue $k^*((\mathbb{Q}))$. 
%Then, we also have $K^* \substruct (K^{**},v^{**})$ (with the same canonical field embedding of $K^*$ to $K^{**}$) since $k^* \substruct k^*((\Qbb))$. We found a field extension $K^{**}$ of $(K,v)$ with  two comparable valuations $u^{**}, v^{**}$ such that $(K,v) \substruct (K^{**},u^{**})$ and $(K,v) \substruct (K^{**},v^{**})$. This concludes our proof.
\end{proof}

%We see that the properties of augmentability by infinitesimals together with the existence of particular power series extensions of fields can be used to apply \Cref{prop:comparablevaluations0Definability} and \Cref{prop:comparablevaluationsDefinability}, and give a characterization of henselian equicharacteristic $0$ valued fields whose valuation is definable (with and without parameters) in the language of rings.  

\subsubsection{Weak augmentability}

\begin{definition}\label{def:WeakAugmentableVF}
    Let $(K,v)$ be a valued field.  We say that $(K,v)$ is \emph{weakly augmentable} if there is a field $K^{*}$ and two valuations $u$ and $w$ on $K^*$ with $\mathcal{O}_u \subsetneq \mathcal{O}_w$ such that $(K,v)\equiv(K^{*},u)\equiv(K^{*},w)$. 
\end{definition}

\begin{observation}\label{obs:maincriterion}
    Rephrasing \Cref{prop:comparablevaluations0Definability}, we have that $(K,v)$ is weakly augmentable if and only if the valuation $v$ is not $\emptyset$-definable in the language of rings $\Lring$.  
    %Indeed, if there is an $\lL_{rings}$-formula $\phi(x)$ without parameters that defines the valuation ring $\Ocal_v$, then it must also define the valuation rings $\Ocal_w$ and $\Ocal_{u}$. Thus, the same formula defines different valuations rings of the same field, which is a contradiction.
\end{observation}

\begin{definition} 
    Let $k$ be a field and $\Gamma$ be an ordered abelian group. We say that $k$ and $\Gamma$ are \emph{weakly co-augmentable} if there %are $k'\equiv k$, $\Gamma'\equiv\Gamma$, and
    is a non-trivial ordered abelian group $\Delta$ such that $k\equiv k((\Delta))$ and $\Gamma\equiv\Gamma\oplus\Delta$.   
\end{definition}

\begin{theorem} \label{thm:augvf}
    Let $(K,k,\Gamma,v)$ be a henselian valued field of equicharacteristic $0$. Then the following are equivalent:
    \begin{enumerate}
    \item\label{item:waugvf1} $v$ is not $\emptyset$-definable
    \item\label{item:waugvf2} $(K,k,\Gamma)$ is weakly augmentable 
    \item\label{item:waugvf3} $k$ and $\Gamma$ are weakly co-augmentable. 
    \end{enumerate}
\end{theorem}
\begin{proof}
The equivalence $(1)\Leftrightarrow(2)$ is essentially \Cref{prop:comparablevaluations0Definability} together with the definition of weakly augmentable valued fields.
We prove $(2)\Leftrightarrow(3)$, starting from the left-to-right implication. Let $K^*$ be a field with two strictly comparable valuations $\mathcal{O}_u\subsetneq \mathcal{O}_w$, such that $$(K,\Gamma,k,v)\equiv(K^*,\Gamma_{u},k_{u},u)\equiv(K^*,\Gamma_{w},k_{w},w)$$
holds. Since $\mathcal{O}_u\subsetneq \mathcal{O}_w$, we have $\Gamma_{w}\cong\Gamma_{u}/\Delta$, for some nontrivial convex subgroup $\Delta$ of $\Gamma_u$.
%Consider an $\aleph_{1}$-saturated extension $$\Delta^{\#}\longrightarrow\Gamma_{u}^{\#}\longrightarrow (\Gamma_{u}/\Delta)^{\#}$$ of the exact sequence $$\Delta\longrightarrow\Gamma_{u}\longrightarrow (\Gamma_{u}/\Delta).$$ Then, $\Gamma_u^{\#} \simeq (\Gamma_u/\Delta)^{\#}\oplus\Delta^{\#}$, since an $\aleph_{1}$-saturated exact sequence of ordered abelian groups splits (\Cref{cor:OagSaturatedExactSequenceSplit}). 
By \Cref{cor:OagSaturatedExactSequenceSplit}, $\Gamma_u \equiv \Gamma_u/\Delta \oplus \Delta$. Therefore, we have $$\Gamma\equiv\Gamma_{u}\equiv\Gamma_{w}\oplus\Delta\equiv\Gamma\oplus\Delta.$$ %i.e. $\Gamma$ is augmentable by $\Delta$.
On the other hand, since $u$ induces a nontrivial henselian valuation on $k_{w}$ with value group $\Delta$ and residue field $k_{u}$, by the Ax-Kochen/Ershov theorem it follows that $k_{w}\equiv k_{u}((\Delta))$. Therefore, we have $$k\equiv k_w\equiv k_u((\Delta))\equiv k((\Delta)),$$ that is, $k$ and $\Gamma$ are weakly co-augmentable.

For the other implication, let $k$ and $\Gamma$ be weakly co-augmentable by a common ordered abelian group $\Delta$. 
Then, by the Ax-Kochen/Ershov principle, $(K,\Gamma, k, v)$ is elementarily equivalent to the field 
$k((\Delta))((\Gamma))$ with two distinct comparable valuations: one being the power series valuation with 
value group $\Gamma \oplus \Delta$ and residue field $k$,
and the other the power series valuation with value group $\Gamma$ and residue field $k((\Delta))$. 
\end{proof}

\begin{remark}
 In \cite{HJ15}, the authors give an example of a valuation ring which is $\emptyset$-definable, but neither with an $\forall\exists$ or $\exists\forall$ formula. This example fits in our framework, since its residue field and value group are not weakly co-augmentable.

 Therefore, the quantifier alternation of the formula defining the valuation in a non weakly co-augmentable henselian valued field can be high.
\end{remark}

%Recalling \Cref{obs:maincriterion}, we make the following remark.
%\begin{remark}
%    Let $(K,v)$ be an equicharacteristic $0$ henselian valued field. Then, the weakly co-augmentability of the residue field and value group is a necessary and sufficient condition for the valuation to be not $\emptyset$-definable in the language of rings.
%\end{remark}

\subsubsection{Strong augmentability}

\begin{definition}\label{def:strongaugmentableVF}
    Let $(K,v)$ be a valued field. We say that $(K,v)$ is \emph{strongly augmentable} if there is a field extension $K^*$ of $K$ and two strictly comparable valuations $\mathcal{O}_u\subsetneq \mathcal{O}_w$ on $K^*$ such that $(K,v)\substruct(K^*,u)$ and $(K,v)\substruct(K^*,w)$.
\end{definition}

\begin{observation}\label{obs:strongmaincriterion}
    Rephrasing \Cref{prop:comparablevaluationsDefinability}, we have that $(K,v)$ is strongly augmentable if and only if the valuation $v$ is not definable (with parameters) in the language of rings $\Lring$.
\end{observation}

\begin{definition}
    Let $k$ be a field and $\Gamma$ be an ordered abelian group. We say that $k$ and $\Gamma$ are \emph{strongly co-augmentable} if %there are elementary extensions $k\substruct k'$, 
    %$\Gamma\substruct\Gamma'$, and 
    there is a non-trivial ordered abelian group $\Delta$ such that $k\substruct k((\Delta))$ and $\Gamma\substruct\Gamma\oplus\Delta$. 
\end{definition}

\begin{remark}\label{rem:soft}\  
\begin{itemize}
    \item Let $k$ be a field of characteristic $0$. If there is an elementary extension $k' \superstruct k$ such that $ k' \substruct k'((\Delta))$, then $k \substruct k((\Delta))$. Indeed, by Ax-Kochen/Ershov, we have
    $k \subseteq k((\Delta)) \substruct k'((\Delta))$,  and thus we also have  $k  \substruct k((\Delta))$. \\
    \item Similarly, if there is an elementary extension $\Gamma' \superstruct \Gamma$ such that $\Gamma' \substruct \Gamma'\oplus \Delta$, then $ \Gamma \subseteq \Gamma \oplus \Delta \substruct \Gamma'\oplus \Delta$, and thus also $\Gamma  \substruct \Gamma \oplus \Delta$.
    %\item When $k$ is of characteristic $0$, one can therefore ``soften'' the definition of co-augmentability as follows:
%$k$ and $\Gamma$ are strongly co-augmentable if and only if there is a non-trivial ordered abelian group $\Delta$ such that $k\substruct k((\Delta))$ and $\Gamma \substruct\Gamma \oplus\Delta$. This is not necessarily true in positive characteristic $p$.
\end{itemize}
\end{remark}

We give the following example of a weakly augmentable valued field which is not strongly augmentable.
\begin{example}\label{ex:VFWeakNotStrongAugmentable}
    Consider the field $k:=\Rbb ((\bigoplus_\omega \Zbb)) \cong \Rbb ((\bigoplus_{1+\omega} \Zbb))$ and the ordered abelian group $\Gamma:=\bigoplus_{\omega^*+1} \Zbb \cong \bigoplus_{\omega^*} \Zbb$. The Hahn series 
    \[K:= k((\Gamma)) := \Rbb ((\bigoplus_\omega \Zbb)) ((\bigoplus_{\omega^*+1} \Zbb))\]
    is weakly augmentable, as it is isomorphic to its coarsening
    \[\Rbb ((\bigoplus_{1+\omega} \Zbb)) ((\bigoplus_{\omega^*} \Zbb)),\]
    but it is not strongly augmentable.
\end{example}
That $\Gamma$ is not strongly augmentable is shown in \Cref{ex:OAGNotStrongAug}. We only note here that the proof essentially relies on the fact that the natural embedding of linear orders $\omega^* \rightarrow \omega^*+1; n \mapsto n$ is not elementary.

\begin{theorem} \label{thm:straugvf}
    Let $(K,k,\Gamma, v)$ be a henselian valued field of equicharacteristic $0$. Then the following are equivalent:
    \begin{enumerate}
    \item\label{item:augvf1} $v$ is not definable
    \item\label{item:augvf2} $(K,k,\Gamma)$ is strongly augmentable 
    \item\label{item:augvf3} $k$ and $\Gamma$ are strongly co-augmentable. 
    \end{enumerate} \end{theorem}
\begin{proof} 
The equivalence $(1)\Leftrightarrow(2)$ is \Cref{prop:comparablevaluationsDefinability} together with the definition of strongly augmentable valued fields. We prove $(2)\Leftrightarrow(3)$, starting from the left-to-right implication. Let $K^*$ be a field extension of $K$ with two valuations $u,w$ with $\mathcal{O}_u\subsetneq \mathcal{O}_w$, such that $$(K,\Gamma,k,v)\substruct(K^*,\Gamma_{u},k_{u},u)$$ and $$(K,\Gamma,k,v)\substruct(K^*,\Gamma_{w},k_{w},w).$$ Since $\mathcal{O}_u\subsetneq \mathcal{O}_w$, then again $\Gamma_{w}\cong\Gamma_{u}/\Delta$, where $\Delta$ is convex in $\Gamma_u$, and $u$ induces a valuation on $k_w$ with value group $\Delta$ and residue field $k_{u}$.

The elementary embedding version of the Ax-Kochen/Ershov theorem implies that $k_{w}\equiv_k k_{u}((\Delta))$. By Keisler-Shelah, there exists a non-principal ultrafilter $\Ucal$ such that $k_{w}$ and $k_{u}((\Delta))$ have isomorphic ultrapowers over $\Ucal$ and the isomorphism is the identity when restricted to $k^\Ucal$, i.e.
$$k_{w}^{\Ucal}\cong_{k^\Ucal} k_{u}((\Delta))^{\Ucal}.$$
By \Cref{fact:EmbedsHahnSeries}, $$k_{u}((\Delta))^{\Ucal}\substruct k_{u}^{\Ucal}((\Delta^{\Ucal})),$$
and since $k^{\Ucal}\substruct k_{u}^{\Ucal}$, we have that $$k^{\Ucal}((\Delta^{\Ucal}))\substruct k_{u}^{\Ucal}((\Delta^{\Ucal})).$$
Then, there is a compatible embedding such that $$k^{\Ucal}\substruct k^{\Ucal}((\Delta^{\Ucal})).$$
On the other hand, consider the exact sequence $$\Delta^{\Ucal}\longrightarrow\Gamma_{u}^{\Ucal}\longrightarrow (\Gamma_{u}/\Delta)^{\Ucal}.$$ Since it is $\aleph_{1}$-saturated, by \Cref{cor:OagSaturatedExactSequenceSplit}, we have $\Gamma_u^{\Ucal} \cong (\Gamma_u/\Delta)^{\Ucal}\oplus\Delta^{\Ucal}.$ It follows that $$\Gamma_{u}^{\Ucal}\cong \Gamma_{w}^{\Ucal}\oplus\Delta^{\Ucal}.$$ Now, by $\Gamma^\Ucal\substruct\Gamma_{w}^\Ucal$, we have that 
$$\Gamma^{\Ucal}\oplus\Delta^{\Ucal}\substruct\Gamma_{w}^{\Ucal}\oplus\Delta^{\Ucal}\cong\Gamma_{u}^{\Ucal}.$$ Moreover, since $\Gamma^\Ucal\substruct\Gamma_{u}^\Ucal$ and $\Gamma^{\Ucal}\subseteq\Gamma^{\Ucal}\oplus\Delta^{\Ucal}$, there is a compatible embedding such that $$\Gamma^{\Ucal}\substruct\Gamma^{\Ucal}\oplus\Delta^{\Ucal}.$$ That is, $k^{\mathcal U}$ and $\Gamma^{\mathcal U}$ are strongly co-augmentable by $\Delta^{\Ucal}$. By \Cref{rem:soft}, $k$ and $\Gamma$ are strongly co-augmentable.

For the other implication, let $k$ and $\Gamma$ be strongly co-augmentable by an ordered abelian group $\Delta$. Then, again, by the elementary embedding version of the Ax-Kochen/Ershov theorem, it follows that $(K,\Gamma, k, v)$ embeds elementarily in the Hahn product $k((\Delta))((\Gamma))$ with the two different comparable valuations: one with residue field $k$ and value group $\Gamma\oplus\Delta$ and the other with residue field $k((\Delta))$ and value group $\Gamma$. 

\end{proof}

%Recalling \Cref{obs:strongmaincriterion}, we make the following remark.
%\begin{remark}
%    Let $(K,v)$ be an equicharacteristic $0$ henselian valued field. Then, the strong co-augmentability of the residue field and value group is a necessary and sufficient condition for the valuation to be not definable (over parameters) in the language of rings.
%\end{remark}

We can develop \Cref{ex:VFWeakNotStrongAugmentable} further to illustrate Theorems \ref{thm:augvf} and \ref{thm:straugvf}:
\begin{example}
    Let $K \coloneqq \mathbb{R}((\bigoplus_\zeta \mathbb Z))= \mathbb{R}((\bigoplus_{\omega} \mathbb Z))((\bigoplus_{\omega^*} \mathbb Z)) $. Then, the only $\emptyset$-definable henselian valuation is $\val : K \rightarrow \bigoplus_\zeta \mathbb Z$. All other henselian valuations are definable over parameters, but not $\emptyset$-definable. This results from the following:
    \begin{itemize}
        \item $\mathbb{R}$ and $\bigoplus_\zeta \mathbb Z$ are not weakly co-augmentable; indeed, $\mathbb R((H))$ is only real-closed if $H$ is divisible, that is, $\mathbb R$ is only augmentable (both weakly and strongly) by divisible ordered abelian groups, while $\bigoplus_\zeta \mathbb Z$ is not (weakly) augmentable by any divisible ordered abelian group.
        \item $\mathbb{R}((\bigoplus_{\omega} \mathbb{Z}))$ and $\bigoplus_{\omega^*} \mathbb Z$ are weakly co-augmentable by $\mathbb{Z}$.
        \item $\mathbb{R}((\bigoplus_{\omega} \mathbb{Z}))$ and $\bigoplus_{\omega^*} \mathbb Z$ are not strongly co-augmentable, since $\bigoplus_{\omega^*} \mathbb Z$ has a  positive element and is therefore not strongly augmentable.
    \end{itemize}
    
    Similarly, we have that any non-trivial convex subgroups of $\bigoplus_\zeta \mathbb Z$ is definable, but not $\emptyset$-definable.
\end{example}
\subsection{Augmentable ordered structures}
To conclude this section, we give two more applications of Beth's definability theorem, this time applied to ordered abelian groups and coloured linear orders. In a very similar way, we can link the definability of convex substructures/end segments with a certain notion of co-augmentability. %We now study the notions of augmentability and co-augmentability in ordered abelian groups and in coloured linear orders. Similarly as for henselian valued fields, these notions relate to definability of certain important sets. 

The results of this subsection are self-contained, in the sense that no other results of the paper rely on them. For the sake of brevity, we decided to omit some parts of the proofs.

\subsubsection{Ordered abelian groups}

\begin{definition}\label{def:CoAugOAGs}
    Let $G$ and $H$ be ordered abelian groups. Then the pair $(G,H)$ is \emph{weakly co-augmentable} if there is a non-trivial ordered abelian group $\Delta$ such that we have:
    \begin{itemize}
        \item $G \equiv G \oplus \Delta$,
        \item $H \equiv \Delta \oplus H$.
    \end{itemize}
    The pair $(G,H)$ is \emph{strongly co-augmentable} if there is a non-trivial ordered abelian group $\Delta$ such that, for the natural embedding, we have
    \begin{itemize}
        \item $G \substruct  G \oplus \Delta$,
        \item $H \substruct \Delta \oplus H$.
    \end{itemize}
\end{definition}

\begin{proposition}\label{prop:co-aug-oag}
    Let $G$ be an ordered abelian group, and $C$ a convex subgroup of $G$. Then 
    
    \begin{itemize}
        \item $C$ is $\emptyset$-definable if and only if $(G/C,C)$ is \emph{not} weakly co-augmentable.
        \item $C$ is definable over parameters if and only if $(G/C,C)$ is \emph{not} strongly co-augmentable.
    \end{itemize}
\end{proposition}
\begin{proof}
We prove the weak case. The proof for the strong case is left to the reader.
To show the right-to-left implication, assume that $C$ is not $\emptyset$-definable. By Beth's definability theorem, there is an elementary extension $G^*$ and convex subgroups $C_1 \trianglelefteq  C_2$ of $G^*$ such that 
\[(G^*,C_1) \equiv (G^*,C_2) \equiv (G,C)\]
in the language $\Loag\cup \{P\}$ with a unary predicate for the convex subgroup. At the cost of taking an ultrapower, we may assume that both $C_1,C_2$ are $\aleph_1$-saturated. Then set $\Delta\coloneq C_2/C_1 $. By \Cref{cor:OagSaturatedExactSequenceSplit}, we have 
\begin{itemize}
    \item $C \equiv C_2 \cong \Delta\oplus C_1 \equiv \Delta\oplus C $,
    \item $G/C \equiv G^*/C_1 \cong G^*/C_2 \oplus \Delta \equiv G/C \oplus \Delta$.
\end{itemize}
In other words, $(G/C,C)$ is co-augmentable by $\Delta$, as wanted. \\
To show the left-to-right implication,  assume that $(G/C,C)$ is co-augmentable by $\Delta$. Let $(G^*,C^*)$ be an $\aleph_1$-saturated extension of $(G,C)$. Then $G^* \cong G^*/C^* \oplus C^*$ and we observe that $(G^*/C^*,C^*)$ is also co-augmentable by $\Delta$:
\begin{itemize}
    \item $G^*/C^* \equiv G^*/C^* \oplus \Delta $,
    \item $C^* \equiv  \Delta \oplus C^* $.
\end{itemize}
Then, we have:
\[ (G^*/C^* \oplus \Delta \oplus C^*, C^*) \equiv (G^*/C^* \oplus C^*, C^*) \equiv (G^*/C^* \oplus \Delta \oplus C^*, \Delta \oplus C^*).\]
A $\emptyset$-$\Loag$-formula cannot define $C$ in $G$, as then it will define both $\Delta \oplus C^*$ and $C^*$  in $G^*/C^* \oplus \Delta \oplus C^*$, which is absurd.
\end{proof}

The following corollary shows that if $(K,v)$ is a henselian valued field of equicharacteristic 0, then the notions of co-augmentability give a one-to-one correspondence between the ($\emptyset$)-definable coarsenings of $v$ in $\Lring$ and the ($\emptyset$)-definable convex subgroups of $vK$ in $\Loag$. %So far, only one implication was known to hold (Pas89), and the proof relied on the value group being stably embedded in the valued field.

\begin{corollary}
Let $(K,v)$ be a henselian valued field of equicharacteristic 0 and $w$ a coarsening of $v$. Let $\Delta_w$ be the corresponding convex subgroup of $vK$. Then $w$ is definable (resp.~$\emptyset$-definable) in $\Lring$ if and only if the pair $(wK,\Delta_w)$ is not strongly (resp.~ weakly) co-augmentable, if and only if $\Delta_w$ is definable (resp.~ $\emptyset$-definable) in $vK$ in $\Loag$.
\end{corollary}
\begin{proof}
    We prove the weak case. The proof of the strong case is left to the reader. By \Cref{prop:co-aug-oag}, it is enough to show the first equivalence. For the left-to-right implication, assume that $(wK,\Delta_{w})$ are co-augmentable by $H$. Recall that, by the Ax-Kochen/Ershov theorem in equicharacteristic $0$, we have $Kw\equiv Kv((\Delta_{w}))$. Then it follows that $$(K,wK,Kw,w)\equiv(Kv((\Delta_{w}))((H))((wK)), wK\oplus H, Kv((\Delta_{w})),u),$$ and $$(K,wK,Kw,w)\equiv(Kv((\Delta_{w}))((H))((wK)), wK, Kv((\Delta_{w}))((H)),u'),$$ where $u$ and $u'$ are two strictly comparable valuations on $Kv((\Delta_{w}))((H))((wK))$. Hence $w$ is not $\emptyset$-definable in the language of rings $\Lring$ by Beth's definability theorem.\\
    For the other implication, assume $w$ is not $\emptyset$-definable in $\Lring$. Then, by \Cref{thm:augvf}, $Kw$ and $wK$ are weakly co-augmentable by a non-trivial ordered abelian group $H$, that is $Kw\equiv Kw((H))$ and $wK\equiv wK\oplus H$. We only need to show that $H$ is an infinite augment of $\Delta_{w}$. Since $Kw\equiv Kv((\Delta_{w}))$, it follows that $$Kw\equiv Kv((\Delta_{w}))((H))\cong Kv((H\oplus\Delta_{w})),$$ that is $\Delta_{w}\equiv H\oplus\Delta_{w}$.
\end{proof}
A comparable statement in the context of almost real closed field can already be found in \cite[Theorem 4.4.]{DF96}. 

\subsubsection{Coloured linear orders}
\begin{definition}
 Let $X$ and $Y$ be coloured linear orders. We say that the pair $(X,Y)$ is weakly (resp. strongly) co-augmentable if there exists a non-empty coloured linear order $Z$ such that both $X\equiv X+Z$ (resp. $X\substruct X+Z$) and $Y\equiv Z+Y$ (resp. $Y\substruct Z+Y$) hold.
\end{definition}

\begin{proposition}\ 
\begin{enumerate}
\item An end segment $I$ of a coloured linear order $X$ is $\emptyset$-definable if and only if the pair $(X\setminus I,I)$ is not weakly co-augmentable.
\item An end segment $I$ of a coloured linear order $X$ is definable with parameters if and only if the pair $(X\setminus I,I)$ is not strongly co-augmentable.
\end{enumerate}
\end{proposition}
The proof is very similar to the previous ones and will be omitted for the sake of brevity.

Working with spines, it is now possible to establish a correspondence between definable convex subgroups of an ordered abelian group $G$ and definable end segments of the coloured linear order $\Sp_n(G)$, as was done by Delon and Farré in \cite[Corollary 4.2]{DF96}. Their result can be rephrased in terms of co-augmentability.

% \begin{corollary}
%  Let $G$ be an ordered abelian group and $C$ a convex subgroup. For $n\geqslant 2$, let $I_n(C)$ be the end segment of $\Sp_n(G)$ defined by $I_n(C)=\sset{x\in\Sp_n(G)}{C<x}$. Then TFAE:
%  \begin{enumerate}
%  \item $C$ is definable (resp.~$\emptyset$-definable) in $G$ in the pure language $\Loag$
%  \item the pair $(G/C,C)$ is not strongly (resp.~weakly) co-augmentable
%  \item for some $n$, we have that $x\in C$ holds if and only if $A_n(x)\in I_n(C)$, and $I_n(C)$ is definable (resp.~$\emptyset$-definable) in $\Sp_n(G)$ in the pure language of coloured linear orders $\Lspn$
%  \item for some $n$, we have that $x\in C$ holds if and only if $A_n(x)\in I_n(C)$, and the pair $(\Sp_n(G)\setminus I_n(C),I_n(C))$ is not strongly (resp.~weakly) co-augmentable.
% \end{enumerate}
% \end{corollary}

%auto-ignore
\section{Augmentability of Ordered Abelian Groups}\label{sec:aug-oag}

We established in the previous section a clear link between definability of valuations and augmentability of ordered abelian groups. It remains to understand exactly when ordered abelian groups are augmentable. As hinted by \Cref{prop:co-aug-oag} and as we will see later, we actually need to consider ``both sides'' of augmentability, that we call \emph{augmentability by infinite} and \emph{by infinitisemals}.
%\footnote{We used ``left'' and ``right'' to distinguished these notions, but this terminology eventually generates some confusion.}:

\begin{definition}\label{def:WeakAugmentsOAGs}
Let $G$ be an ordered abelian group. We say that
\begin{itemize}
    \item $G$ is \emph{weakly-augmentable by infinitesimals} if  $G\equiv G\oplus H$ for some non-trivial ordered abelian group $H$.
    \item $G$ is \emph{weakly-augmentable by infinites} if  $G\equiv H \oplus G$ for some non-trivial ordered abelian group $H$.
\end{itemize}
\end{definition}

\begin{definition}\label{def:StrongAugmentsOAGs}
Let $G$ be an ordered abelian group. We say that
\begin{itemize}
    \item $G$ is \emph{strongly-augmentable by infinitesimals} if $G\substruct G\oplus H$ for some non-trivial ordered abelian group $H$.
    \item $G$ is \emph{strongly-augmentable by infinites} if $G\substruct H\oplus G$ for some non-trivial ordered abelian group $H$.
\end{itemize}
\end{definition}

We call such $H$ strong- (resp.~weak-)augments of $G$. We will mainly rely on the work of Delon and Lucas \cite{DelonLucas} to characterize strongly-augmentable ordered abelian groups. We will need, however, to delve deeper into the theory of spines to characterize the class of ordered abelian groups weakly-augmentable by infinitesimals. For that, we will use \emph{Schmitt's Habilitationsschrift} \cite{Sch82} as a main reference.

%auto-ignore
\subsection{Strong Augmentability by infinites}

We show that all ordered abelian groups are strongly augmentable by infinites. This can be easily deduced by the first Proposition of \cite{DelonLucas}:

\begin{fact}[See {\cite[Proposition 1]{DelonLucas}}]\label{lem:ConvxHull}
    Let $G$ be an ordered abelian group and $G'$ an elementary extension. Then the convex hull $\CVX_{G'}(G)$ of $G$ in $G'$ is also an elementary extension of $G$ and an elementary substructure of $G'$:
    \[ G \substruct \CVX_{G'}(G) \substruct G'.\]
\end{fact}

\begin{remark}\label{rmk:AI-is-elementary}
Let $G,H,H'$ be ordered abelian groups. If $G\substruct H\oplus G$ and $H'\equiv H$, then $G\substruct H'\oplus G$, as $H\oplus G\equiv_G H'\oplus G$.

Similarly, if $G\substruct G\oplus H$ and $H'\equiv H$, then $G\substruct G\oplus H'$.
\end{remark}

\begin{theorem}\label{thm:strongleftaugmentabilityOAG}
    All non-trivial ordered abelian groups are augmentable by infinites.
\end{theorem}

\begin{proof}
Let $G$ be a non-trivial ordered abelian group and let $G'$ be an elementary extension of $G$ realising a type at $+\infty$. Let $H$ be the convex hull of $G$ in $G'$. Then $G\substruct H$ by \Cref{lem:ConvxHull} and $G'/H\neq\{0\}$. Consider the exact sequence $0\rightarrow H\rightarrow G'\rightarrow G'/H\rightarrow 0$. Move if needed to an $\aleph_1$-saturated extension $(G^*,H^*)\superstruct (G',H)$. By \Cref{cor:OagSaturatedExactSequenceSplit}, $G^*\cong G^*/H^*\oplus H^*$. Since $H\substruct H^*$, then also $G\substruct H^*$. This implies that $A\oplus G\substruct A\oplus H^*$ for any ordered abelian group $A$, in particular, $G^*/H^*\oplus G\substruct G^*/H^*\oplus H^*\cong G^*$. Now since $G\substruct G^*$ and $G\subseteq G^*/H^*\oplus G$, we have that $G\substruct G^*/H^*\oplus G$.
\end{proof}

Note that it follows from the proof that all ordered abelian groups $G$ admit a
proper elementary extension $G^*$ such that $G$ is a convex subgroup of $G^*$.

We can, in fact, be more precise and characterize when an ordered abelian group admits an $n$-divisible infinite augment. The statement (and its proof) use the formalism of spines, which we will introduce immediately after.

\begin{corollary}
   Let $G$ be an ordered abelian group. Then for any natural number $n\geqslant 2$, TFAE:
   \begin{enumerate}
   \item Any infinite strong-augment of $G$ is $n$-divisible,
    \item Some infinite strong-augment of $G$ is $n$-divisible,
   \item $\Sp_n(G)$ has an initial point.
   \end{enumerate}

   Furthermore, TFAE:
   \begin{enumerate}\setcounter{enumi}{3}
   \item $\mathbb Q$ is an infinite strong-augment of $G$,
   \item Any divisible ordered abelian group is an infinite strong-augment of $G$,
   \item Any infinite strong augment of $G$ is divisible,
   \item For any $n$, $\Sp_n(G)$ has an initial point.
   \end{enumerate}
   \label{cor_div}
\end{corollary}

\begin{proof}
The implication $(1)\Rightarrow(2)$ is immediate. We prove $(3)\Rightarrow(1)$. Let $H$ be a non $n$-divisible ordered abelian group and assume $G\substruct H\oplus G$. Let $h\in H\notin nH$, then in $H\oplus G$, $F_n(h,0)\supseteq G$. Recall that $\Sp_n(G)\substruct \Sp_n(H\oplus G)$ (see \cite[Lemma 0]{DelonLucas}); thus, \[\Sp_n(G)\vDash (\forall x)(\exists y) (y<x),\]
that is, $\Sp_n(G)$ doesn't have an initial point.\\

To show $(2)\Rightarrow(3)$, assume now that $H$ is $n$-divisible and that $G\substruct H\oplus G$. %Taking any $h\in H\setminus\{0\}$ and any $g\in G$, we have that $A_n((h,g))\subseteq G$. 
Since $H$ is $n$-divisible group and therefore $n$-regular, all elements $(h,g)$ and $(h',g')$, with $h,h'\in H\setminus\{0\}$ and $g,g'\in G$, belong to the same $n$-regular class. This mean that $A_n((h,g))=A_n((h',g'))$ for all $h,h'\in H\setminus\{0\}$ and $g,g'\in G$. Finally, $F_n((h,g))\subseteq A_n((h,g))$ by definition. Therefore, $A_n((h,g))$ for $h\in H\setminus\{0\}$ is an initial point in the $n$-spine of $H\oplus G$, and by elementarity such an initial point also exists in $\Sp_n(G)$.

\ 

The second set of equivalences follows immediately from the first and from \Cref{rmk:AI-is-elementary} -- since the $\Loag$-theory DOAG of divisible non-trivial ordered abelian groups is complete.
\end{proof}

Note that this results characterizes infinite strong augments, but says nothing about weak augments. For example, the $n$-spine of $\bigoplus_{\omega}\mathbb Z$ has an initial point for every $n$, and thus, any infinite strong augment is divisible, by the previous result. However, $\mathbb Z$ is also an infinite weak augment of $\bigoplus_{\omega}\mathbb Z$.

%auto-ignore
\subsection{Review of the language of Schmitt}\label{subsec:lang-of-spines}
    We present here the language of spines, as developed by Schmitt in his \emph{Habilitationsschrift} \cite{Sch82}. We were also guided by the excellent exposition made in \cite{DelonLucas}. Notice that (minor) choices of convention may be different to the literature.
    
    Let $G$ be an ordered abelian group. For $g\in G$, we define the archimedean class of $g$ as follows:
    \begin{itemize}
        \item $A(g,G)$ denotes the largest convex subgroup of $G$  not containing $g$, 
        \item $B(g,G)$ denotes the smallest convex subgroup of $G$ containing $g$,
        \item $R(g,G)\coloneq B(g,G)/A(g,G)$ the archimedean class of $g$ in $G$.
    \end{itemize}
    Next, we define the $n$-regular class of $g$:
    \begin{definition}\label{def:spines}
        Let $n\geqslant 2$ and $g\in G$. \begin{itemize}
            \item We let $A_n(g,G)$ denote the convex subgroup:
        \[ \bigcap \{C \trianglelefteq G : B(g,G)/ C \text{ is $n$-regular}\}\]
            \item We let $B_n(g,G)$ denote the convex subgroup:
        \[ \bigcup \{C \trianglelefteq G : C/A(g,G) \text{ is $n$-regular}\}.\] 
            \item We let $R_n(g,G)$ denote the quotient $B_n(g,G)/A_n(g,G)$. This is an $n$-regular group, and in fact the largest one containing a projection of $g$.
        \end{itemize}
    To an element $g$, we can associate another important convex subgroup, called the $n$\emph{-fundament of} $g$, denoted $F_n(g,G)$:
        \[ F_n(g,G) \coloneqq \bigcup \{C \trianglelefteq G : C \cap (g + nG)= \emptyset\}\]    
     If $g\notin nG$, then $F_n(g,G)$ is the largest convex subgroup $C$ such that $g+C$ is not $n$-divisible in $G/C$. If $g\in nG$, we sometimes write $F_n(g,G) = \emptyset$.
     
    \end{definition}
       
    \begin{definition}\label{def:Spn1}
            The $n$-spine $\Sp_n(G)$ of $G$ is the set of definable convex subgroups 
            \[\{A_n(g,G) : g \in G, g\neq 0 \} \cup \{F_n(g,G) : g \in G \setminus nG \}\footnote{Notice that $F_n(g,G)= \emptyset$ for all $g\in nG$, but $\emptyset$ is not considered as an element of the spine by our definition. We deviate from Schmitt's convention in \cite{Sch82} here for a good reason: otherwise, we would have a last element in all the spines, which would make a characterization of augmentability more difficult to state.}  \]
            with the reverse ordering for the inclusion: if $C,C' \in \Sp_n(G)$,
            \[ C \leqslant C' \Leftrightarrow \ C \supseteq C'.\]
            It will also be equipped with colours that we define below.
    \end{definition}    
 If there is no ambiguity on the ambient group, we may simply write $A_n(g)$, $F_n(g)$, $B_n(g)$, etc.

    \begin{definition}\label{def:alphas}
        For $C$ a convex subgroup of $G$, we set
    \[\Gamma_1(n,C,G) \coloneq \{g \in G : F_n(g,G) \subseteq C \}\]
    \[\Gamma_2(n,C,G) \coloneq \{g \in G : F_n(g,G) \subsetneq C \}\]
    and 
    \[\Gamma(n,C,G) \coloneq \Gamma_1(n,C,G)/\Gamma_2(n,C,G)\]

    We fix the group $G$ and omit it in the notation, by writing $\Gamma(n,C),\Gamma_1(n,C)$ and $\Gamma_2(n,C)$.  
    For $\Gamma$ an abelian group, $p$ a prime and $s$ a natural number, we denote by $\Gamma[p]$ the $p$-torsion of $\Gamma$, and we define 
    \[\alpha_{p,s} (\Gamma) \coloneq \begin{cases}\dim_{\mathbb{F}_p} ((p^{s-1}\Gamma)[p]/(p^s\Gamma)[p]) & \text{ if finite} \\ \infty & \text{otherwise} \end{cases}.\]     This is the number of copies of the cyclic group $C_{p^s
    }$ in $\Gamma$.
    \end{definition}
    
%%% DO NOT ERASE %%%%%%%
%    \begin{fact}[{\cite[Theorem 1.1]{Gur64}}]\label{fact:Gur1.1}
%    We have \[\alpha_{p^s}(\Gamma(p^k,C))= \alpha_{p^s}(\Gamma(p^{k+1},C))\] for $s<k$, and 
%        \[\tag{$*_\alpha$}\alpha_{p^k}(\Gamma(p^k,C))= \alpha_{p^k}(\Gamma(p^{k+1},C)) + \alpha_{p^{k+1}}(\Gamma(p^{k+1},C)) \]
%    for every $k$. 
%    \end{fact}

    \begin{definition}\label{def:Spn2}

    If $R$ is an $n$-regular ordered abelian group and $p|n $, we denote  
    \[\beta_p (R) \coloneq \begin{cases}\dim_{\mathbb{F}_p} (R/p R) & \text{ if finite} \\ \infty & \text{otherwise} \end{cases}.\]
    This is the number of copies of $\mathbb{Z}_{(p)}$ in $R$.
    \end{definition}

    \begin{definition}\label{def:Spn3}
        Let $n \in \mathbb{N}$ and $p^{k_p}$ be largest prime power of $p$ factor of $n$. We equip $Sp_n(G)$ with the following unary predicates $F_n,A_n,D,
        \alpha_{p,s,m},\beta_{p,m}$ , where $s<k_p$ and $m \in \mathbb{N}$, where: 
        \begin{itemize}
            \item $A_n (C)$ if $C= A_n(g)$ for some $g \in G$,
            \item $F_n (C)$ if $C= F_n(g)$ for some $g \in G$,
            \item $D(C)$ if $C=A_n(g)$ for some $g\in G$ and $R_n(g)$ is discrete,
            \item $\alpha_{p,s,m} (C)$ if $C=F_n(g)$ for some $g\in G$ and $\alpha_{p,s} (\Gamma(p^{k_p},C))\geq m$,
            \item $\beta_{p,m}(C)$ if $C=A_n(g)$ for some $g\in G$  and for $g$ such that  $\beta_{p} (R_p(g))\geq m.$
            
        \end{itemize}
        We denote by $\mathcal{L}_{\Sp_n}$ the following language
        \[\{<, A_n,F_n\} \cup \{ \alpha_{p,s,m}: p\vert n, s<k_p, \ m \in \Nbb  \}  \cup \{\beta_{p,m} : p \vert n, \ m\in \Nbb \}.\] 
    
    \end{definition}

Interested readers can find some examples described in \cite{DelonLucas} and \cite{Sch84}.  
We want to emphasise some more. The following example show that $A_n$ and $F_n$ don't need to intersect:
\begin{remark}
        Alternative definitions of the $\alpha_{p^s}(\Gamma(p^{k_p},C))$ can be found in \cite{CH11}, where they are denoted $d_{p,s}(C)$ and $d_{p,s}^\infty(C)$.
    \end{remark}

    \begin{example}
        Recall that $\eta$ denotes the dense linear order $(\mathbb{Q},<)$. Consider $H \coloneq \bigoplus_\eta \Qbb$ and for every irrational $\pi$, a sequence $(r^\pi_i)$ of rational greater than $\pi$ converging to $\pi$. We consider the element $a^\pi =(a^\pi_r)_{r\in \eta} \in \hprod_\eta \Qbb$ where:
        \[a^\pi_r \coloneq 
        \begin{cases}
            1 & \text{ if $r=r^\pi_i$ for some $i$,}\\
            0 & \text{ otherwise.}
        \end{cases}\]
        Set $G\coloneq \braket{H,a^\pi : \pi\in \Rbb\setminus \Qbb }.$ Then, for every integer $n$, $\Sp_n(G)= \Rbb$, $A_n(G) = \Qbb$, and $F_n(G) = \Rbb \setminus \Qbb$. For $x\in A_n(G)$, $R_x = \mathbb{Q}$ and for $x\in F_n(G)$, $\alpha_{p,s}(x)=0$ for $s<k-1$ and $\alpha_{p,k-1}(x)=1$.
        
    \end{example}

    The next example builds on the previous one, hinting that the predicates $\alpha_{p,s}$ can be intricate.  
    \begin{example}
        Let $G$ be as above. Fix an integer $s$, and let $(t_n)$ be an increasing sequence of irrational numbers converging to $0$. Set $b= \sum_{n\in \Nbb} p^sa^{t_n}$. We can insure that $b \in \hprod_\eta \Qbb $ by forcing the support of $a_{t_{n+1}}$ to be in $]t_n,t_{n+1}[$. Now set $G' \coloneqq \braket{G,b}$. Then $F_n(G')$ has an extra point for $x \coloneqq F_n(b)=\CVX_{G'}(\sum_{\geq 0}\mathbb Q)$, and $\alpha_{p,s}(x)=1$, $\alpha_{p,r}(x)=0$ for all $r\neq s$.
    \end{example}
      \begin{fact}[Schmitt]\label{prop:SpinesInterpretableInOAG}\ 
        \begin{itemize}
            \item For all $n\in \mathbb{N}_{>1}$, the $n$-spine $\Sp_n(G)$, equipped with its colors and ordering, is interpretable by $G$ in the ordered group language $\Loag$.
            \item \cite[Corollary 4.7]{Sch82} Let $G_1$ and $G_2$ be two ordered abelian groups. Then $G_1\equiv G_2$ if and only if for all $n\in \mathbb{N}_{>1}$:
            \[\Sp_n(G_1) \equiv \Sp_n(G_2)\]
        \end{itemize}
        
    \end{fact}
     Schmitt proves, in fact, a quantifier elimination result relative to the spine \cite[Theorem 1.7]{Sch84}. We will not state it, as we will use instead the work of Delon and Lucas \cite{DelonLucas}, which is based on this stronger result. 
     
     %Interpretability of the $n$-spine is proved in Chapter 3 of Schmitt's Habilitationschrift. %We will only recall the proof of 
     In the next lemma, we gather some intermediate results from Schmitt's Habilitationschrift that we will use later:
\begin{lemma}\label{lem:FnSChmitt}\ 
Let $n>1$ be an integer, and $G$ be an ordered abelian group.
    \begin{enumerate}
        \item  for $a\in G\setminus nG$,  we have the following equivalence:
    \[g\in F_n(a) \Leftrightarrow (\forall h \in G) (\vert h\vert <n \vert g \vert \rightarrow a+h \notin nG).\]
        \item $F_n(g)= \bigcup_{p^m \vert n}F_{p^m}(g)$,
        \item for all $m\vert n$ such that $n\vert m^k$ for some $k$, and for all $g\in G\setminus nG$, there is $h$ such that $F_n(g)=F_m(h)$.
    \end{enumerate}
\end{lemma}
See \cite[Lemma 2.10]{Sch82} for (2) and (3). 
A proof of first-ordered definability of $F_n$ is scattered in \cite[Chapter 3]{Sch82}, or \cite[p. 393]{Sch84}. However, the exact statement (1) is \cite[Lemma 2.1]{CH11}. %We rewrote a proof of (1) to match the notation of \cite{Sch82}:

%\begin{proof}
%    By definition, $F_n(a)$ is the largest subgroup $C$ such that $a+C$ is not $n$-divisible. In particular,  ``$\Rightarrow$'' is immediate. To show ``$\Leftarrow$'', let $g\in G$ and assume that 
%    \[(\forall h \in G) (\vert h \vert<n \vert g\vert \rightarrow a+h \notin nG)\]
%    holds. Then we need to show that, if $C=\CVX\langle g\rangle$ is the convex subgroup generated by $g$, then $a+C$ is not divisible by $n$. If it were, let $b$ be such that $\delta\coloneqq a+nb \in C$. Without loss of generality,  $g>0$ and $\delta>0$ and there is an integer $m$ such that $(m-1)g<\delta<mg$. Write $m=qn+r$ with $0\leq r< n$. We get  
%    \[0<\delta-qng<ng,\]
%   which is a contradiction.
%\end{proof}

\subsection{Strong Augmentability by infinitesimals}
%Let $G$ be an ordered abelian group. Recall that if $h\in G$ is an element, we let $\vert h \vert $ denote $\max \{h,-h\}$.
To characterize ordered abelian groups which are strong augmentable by infinitesimals, we rely on the second proposition of Delon-Lucas:
\begin{fact}[{\cite[Proposition 2]{DelonLucas}}] \label{fact:DelonLucas}
    Let $H\preccurlyeq G$ be a pair of ordered abelian groups. Let $H_0$ be the $H$-infinitesimals of $G$:
    \[H_0 \coloneq \{g\in G : \forall h \in H (h\neq 0 \rightarrow \vert g\vert < \vert h \vert ) \}.\]
    Assume that for all integers $n$:
    \begin{itemize}
        \item for all $g\in G$, $F_n(g,G)\neq \{0\},$
        \item for all $g\in G$, $F_n(g,G)\neq H_0.$
    \end{itemize}
    Then $H \preccurlyeq G/H_0$ with the obvious embedding.
\end{fact}

\begin{corollary}
     Assume that for all $g\in G$ and all $n$, $F_n(g,G)\neq \{0\}$. Then  $G$ is strongly augmentable by infinitesimals. 
\end{corollary}
\begin{proof}
    Observe that in particular, $G$ has no smallest positive element, and we can find an elementary extension  $G^1$ of $G$ with a non-trivial convex subgroup of infinitesimals $G_0^1$. %By the second condition, there is $n$ such that $G_0^1$ is not divisible by $n$ and in particular, $Sp_n(G_0^1)$ is not empty.  
    \begin{claim}\label{claim}
        There is an elementary extension $G^*$ of $G$ such that for all $n\in \mathbb{N}$ and $g\in G^*$, $G^*_0 \neq F_n(g,G^*)$. 
    \end{claim}
    \begin{proof}[Proof of \Cref{claim}]
        
    If for all $g\in G^1$ and all $n$, $F_n(g,G)\neq G_0^1$, we are done. Assume otherwise and for example that $F_{n_1}(a_1,G^1)= G_0^1$ for some $a_1\in G^1$ and $n_1\in \mathbb{N}$.
    Consider the type over $G \cup \{a^1\}$:  
    \[\rho_1(x)\coloneq \{x \notin F_{n_1}(a_1,G^1)\} \cup \{ \vert x\vert< \vert g\vert : g \in G\setminus\{0\}\}.\]
    This type says that $x$ is a new infinitesimal of $G$, larger than $G_0^1$.
    Let $G^2$ be an extension realising $\rho_1$. Let $G_0^2$ be the subgroup of infinitesimals in $G^2$. Note that $F_{n_1}(a_1,G^2)\neq G_0^2$ since $x\in G_0^2\setminus F_{n_1}(a_1,G^2)$.
    
    If $G_0^2$ is not an $F_n$, then we are done. Else, there is $a_2\in G^2$ and $n_2\in\mathbb N$ such that $F_{n_2}(a_2,G^2)=G^2_0$, and we define similarly 
\[\rho_2(x)\coloneq \{x \notin F_{n_2}(a_2,G^2)\} \cup \{ \vert x\vert< \vert g\vert : g \in G\setminus\{0\}\}.\]
 and $G^3$ realising it. We proceed by induction and define if needed $G^m,a_m,n_m$ and $\rho_m$ for any $m<\omega$.
    If this process stops at some finite point, we are done. Otherwise, let $G^*$ denote the union $\bigcup_m G^m$:
        \begin{center}
        \begin{tikzpicture}
        \draw[thick] (-2,0) -- (1.5,0);
        \draw (0,0) node[anchor=north]{$G$} ;
        \draw (2.2,0) node{$\dots$} ;
        \draw[thick] (2.5,0) -- (3,0);
        \draw (2.5,-0.4) --++ (0,0.8);
        \draw (2.5,0.5) node[anchor=south]{$G_0^3$} ;
        \draw [decorate,
    decoration = {brace}] (7,-0.4) -- ++ (-5,0);
    \draw (4.3,-0.5) node[anchor=north]{$G_0^*$} ;
    \filldraw (7,0) node[anchor=south]{$0$} circle(1pt);

        \draw[thick] (3.3,0) -- (4.3,0);
        \draw (3.3,-0.4) --++ (0,0.8);
        \draw (3.3,0.5) node[anchor=south]{$G_0^2$} ;
        \draw[thick] (4.6,0) -- (6.6,0);
        \draw (4.6,-0.4) --++ (0,0.8);
        \draw (4.6,0.5) node[anchor=south]{$G_0^1$} ;
    \end{tikzpicture}
    \end{center}

    We claim that for all element $g$ in $G^*$ and $n \in \mathbb{N}$, $F_n(g,G^*) \neq G^*_0$.

    Indeed, take $g\in G^*$ and $n\in\mathbb N$. There is $m$ such that $g\in G^m$. Now either 
    \begin{itemize}
        \item $F_n(g,G^m)\supsetneq G^m_0$, which means by definition that it contains an element of our original group $G$. This also holds in $G^*$ since $G^m\substruct G^*$. This means that $F_n(g,G^*)\supsetneq G^*_0$. Or,
        \item $F_{n}(g,G^m)\subseteq G^m_0$. But $G^m_0=F_{n_m}(a_m,G^m)$, and the first-order sentence ``$F_{n}(g,G^m)\subseteq F_{n_m}(a_m,G^m)$'' is true in $G^m$. Moving to $G^{m+1}$, there is an $x$ realising $\rho_m$, which means $x \in G^{m+1}_0 \setminus F_{n_m}(a_m,G^{m+1})$. Since $F_{n_m}(a_m,G^{m+1})\supseteq F_n(g,G^{m+1})$, we have in particular $x\notin F_n(g,G^{m+1})$. Consequently, $x\in G^*_0\setminus F_n(g,G^*)$, that is, $F_n(g,G^*)\neq G^*_0$.
    \end{itemize}

    \renewcommand{\qedsymbol}{$\square_{\text{of \Cref{claim}}}$}
    \end{proof}

    By \Cref{fact:DelonLucas}, $G \substruct G^*/G^*_0$ and in particular
    \begin{align}
        G + G^*_0 \preccurlyeq G^*/G^*_0 \oplus G^*_0. \label{eq:G+H0}
    \end{align}

    Consider the pair $(G^*,G)$, in the language of groups with a predicate for $G$. In this language, $G_0^*$ is definable and the statement ``$G \substruct G^*/G^*_0$'' is therefore  part of the theory of the pair $(G^*,G)$.
    We can therefore move, if needed, to an elementary extension of the pair $(G^*,G)$, and assume that $G^*$,$G_0^*$ and $G^*/G_0^*$ are $\aleph_1$-saturated and \begin{align}
        G \preccurlyeq G^*\cong G^*/G_0^* \oplus G_0^*.\label{eq:G}
    \end{align}
           From \Cref{eq:G,eq:G+H0}, it follows that 
    \[G \preccurlyeq G \oplus G_0^*,\]
    and this shows that $G$ is augmentable by infinitesimals.
    
\end{proof}

As we shall see now, it is, in fact, an equivalence. We can also reformulate the assumption as a property of the spines. We have:

\begin{theorem}\label{thm:StronglyAugmentableInfinitesimalsOAGs}
Let $G$ be an ordered abelian group. The following are equivalent:
\begin{enumerate}
    \item $G$ is strongly augmentable by infinitesimals 
    \item for all $g\in G$ and all $n\geqslant2$, $F_n(g,G)\neq \{0\}$
    \item for all $n\geqslant2$, $\Sp_n(G)$ \emph{doesn't} have a last element satisfying $F_n$.
\end{enumerate}
\end{theorem}

\begin{proof}
 $(2) \Rightarrow (1)$ is the previous result. For $(1) \Rightarrow (2)$, assume that there is $H\neq\{0\}$ such that $G\substruct G\oplus H$ and fix $n\in\mathbb N$ and $g\in G$. If $F_n(g,G)=\{0\}$, then $F_n((g,0),G\oplus H)$ is also $\{0\}$. By definition of $F_n$, this means that $(g,0)$ is $n$-divisible modulo any non-trivial convex subgroup. In particular, $g$ is $n$-divisible in $G=(G\oplus H)/H$; but this means $F_n(g)=\emptyset$ which contradicts the assumption that $F_n(g)=\{0\}$.
 Finally, $(2) \Leftrightarrow (3)$ is immediate: if $\Sp_n(G)$ has a smallest element satisfying $F_n$, then this must correspond to the convex subgroup $\{0\}$, and there is therefore $g\in G$ such that $F_n(g,G)=\{0\}$.  
\end{proof}

This completely characterizes ordered abelian groups which are strongly augmentable by infinitesimals.

\begin{example}\label{ex:OAGNotStrongAug}
    Consider $G\coloneq \bigoplus_{\omega^*} \mathbb{Z}$. Then for all $n\geqslant2$, the order type of $\Sp_n(G)$ is $\omega^*$. All points satisfy exactly the following predicates: 
    \[\forall x \in \Sp_n(G) ~ F_n(x) \wedge A_n(x) \wedge D(x)\wedge \beta_{p,1}(x)\wedge \alpha_{p,v_p(n)-1,1}(x)\]
    where $v_p(n)$ is $p^{th}$ valuation of $n$.

    In particular, $\Sp_n(G)$ has a last element satisfying $F_n$, and $G$ is not strongly augmentable by infinitesimal by \Cref{thm:StronglyAugmentableInfinitesimalsOAGs}. Notice however that $\mathbb Z$  is a weak infinitesimal augment of $G$: we have, indeed, that
    $G \equiv G \oplus \mathbb{Z}$, and in fact, that $G \cong G \oplus \mathbb{Z}$.
\end{example}

Like in the case of 
augmentation by infinites, we 
now characterize ordered abelian groups strongly augmentable by a divisible infinitesimal augment. Parts of this characterization already appear in \cite[Lemma 5.9]{FJ15}.

\begin{lemma}\label{DOAG-aug}
Let $G$ be an ordered abelian group. TFAE:
\begin{enumerate}
\item\label{1} $G\substruct G\oplus H$ for some non-trivial divisible ordered abelian group $H$.
\item\label{2} $G\substruct G\oplus H$ for any divisible ordered abelian group $H$.
\item\label{3} Given an ordered abelian group $H$, we have $G\substruct G\oplus H$ if and only if $H$ is divisible.
\item\label{4} $G\equiv G\oplus H$ for some non-trivial divisible ordered abelian group $H$.
\item\label{5} $G\equiv G\oplus H$ for any divisible ordered abelian group $H$.
\item\label{8} $G$ admits for each prime $p$ a non-trivial $p$-divisible convex subgroup.
\item\label{7} For each prime $p$, there is $g\in G$ such that the interval $(0,g)$ is non-empty and $p$-divisible.
\item\label{9} There is $G'\equiv G$ admitting a non-trivial divisible convex subgroup.
\item\label{0} For each $n\geqslant 2$, $\Sp_n(G)$ has a final point satisfying $\neg F_n$. 
\end{enumerate}
\end{lemma}

\newcommand{\p}[1]{(\ref{#1})}

\begin{proof}
$\p2\Rightarrow\p1$ is clear and $\p1\Rightarrow\p2$ follows from \Cref{rmk:AI-is-elementary} and the completeness of the theory of divisible ordered abelian groups.

Similarly, we get $\p4\Leftrightarrow\p5$.

We now prove $\p8\Rightarrow\p7 \Rightarrow\p9 \Rightarrow\p0 \Rightarrow\p8$. First, assume that \p8 holds. Then $G$ admits a non-trivial $p$-divisible convex subgroup, and any $g>0$ in this convex subgroup fulfills \p7. Next, if \p7 holds, then $G$ has an infinitesimal type, in the sense that $\sset{0<x<\vert g\vert}{g\in G\setminus\{0\}}$ is a finitely satisfiable type -- this is because a $p$-divisible interval is dense. Take $G'$ realizing this type and consider the subgroup $G'_0$ of $G$-infinitesimals in $G'$. It is non-trivial, and it is divisible, as for $p$ and $g\in G$ given as is \p7, $G'$ satisfies ``$0<x<g\rightarrow x$ is $p$-divisible'', that is, we have \p9.

Continuing from \p9 towards \p0, assume \p9 holds and take $G'\equiv G$ with a non-trivial divisible convex subgroup, call it $C$. Fix $n\geqslant 2$. Now any $c\in C$ satisfies $A_n(c)=\{0\}$. This is a final point in $\Sp_n(G')$. Take now $g\in G'$ and compute $F_n(g)$: if $g\in nG'$, then $F_n(g)=\emptyset$ is not in the spine, and if $g\notin nG$, then also $g+c\notin nG$ for any $c\in C$ since $c\in C$ is divisible, that is, $F_n(g)\supseteq C$. In both cases, $F_n(g)\neq\{0\}$, so $\Sp_n(G')$ satisfies ``I have a final point which is $\neg F$''. Since $\Sp_n(G')\equiv\Sp_n(G)$, \p0 holds.

Finally, assume \p0. Fix $p$ prime, by \p0, $\Sp_p(G)$ has a final point which is not an $F_p$. Then this final point must be an $A_p$ %.If $A_p(g)\supsetneq\{0\}$ for some $g\in G$, then taking $h\in A_p(g)\setminus\{0\}$, we have $A_p(g)\supsetneq A_p(h)$, so a final point of colour $A$ in the spine has to correspond to
and for some $g\in G$, we have $A_p(g)=\{0\}$. Now $B_p(g)$ is a convex subgroup of $G$. Take $h\in B_p(g)$, then $F_p(h)\subseteq A_p(h)=A_p(g)=\{0\}$ by properties of the spine. Since the final point is not an $F_p$, we can't have $F_p(h)=\{0\}$; thus, $F_p(h)=\emptyset$, that is, $h$ is $p$-divisible, and $B_p(g)$ is a non-trivial convex subgroup of $G$ which is $p$-divisible. This gives $\p0\Rightarrow\p8$.

To close the loop, observe that $\p3\Rightarrow\p1\Rightarrow\p4\Rightarrow\p9$. Now assume (8) holds, i.e. there is $G'\equiv G$ admitting a non-trivial divisible convex subgroup, name this subgroup $C$. Then $G'\equiv G'/C\oplus C$ by \cref{cor:OagSaturatedExactSequenceSplit}.
To show the right-to-left implication of \p3, take $H$ a divisible ordered abelian group. By properties of divisible ordered abelian groups, $C\substruct C\oplus H$. Thus $G'/C\oplus C\substruct G'/C\oplus C\oplus H$, and since $G\equiv G'/C\oplus C$, we have $G\substruct G\oplus H$.

We show now the right-to-left implication of \p3. Assume that $G\substruct G\oplus H$. Now also $G'\substruct G'\oplus H$. Take $c\in C$, now $G'$ satisfies that any element in the interval $(0,c)$ is divisible, so $G'\oplus H$ also, and in particular, any element in $H$ is divisible.

Thus, we have $\p9\Rightarrow\p3$, which completes the proof.
\end{proof}

Note that this characterizes groups which are both weakly and strongly augmentable by divisible infinitesimals; however, a ``weak'' version of \p3, namely 
\begin{itemize}
    \item[(3')] Given an ordered abelian group $H$, we have $G \equiv G \oplus H$ if and only if $H$ divisible.
\end{itemize}
is stronger than all the above statements. Indeed, we have the following example:
\begin{example}\label{ex:aug-by-Q-and-Z+Q}
$\bigoplus_{\omega^*} (\mathbb Z\oplus\mathbb Q) \equiv (\bigoplus_{\omega^*}\mathbb Z) \oplus\mathbb Q $ is weakly augmentable by $\mathbb Q$. It is, however, also weakly augmentable by $\mathbb Z\oplus\mathbb Q$.
\end{example}

\subsection{Axiomatisation of the theory of spines}\label{subsec:theory-of-spines}
        
We will proceed in the next subsection to characterize ordered abelian groups which are weakly augmentable by infinitesimals. We need to recall and discuss Schmitt's axiomatisation of $n$-spines. We fix an integer $n$ and for each prime $p|n$, we let $p^{k_p}$ be the largest prime power dividing $n$. We sometimes drop the index $n$ in our notation for more visibility, and in particular write $F(x) $ and $ A(x)$ instead of $F_n(x)$ and $A_n(x)$.
        
        \begin{definition}
            Let $T_{\Sp,n}$ be the following list of axioms in the language $\mathcal{L}_{\Sp,n}:$
            \begin{enumerate}
                \item[(AS1)] $<$ is a linear order
                \item[(AS2)] \label{Axiom:Sp2} $(\forall x)~ (\neg\beta_{p,m}(x)\rightarrow \neg \beta_{p,m+1}(x))$ \quad for all $p|n$ and $m\geq 1 $
                \item[(AS3)]  \label{Axiom:Sp3} $(\forall x)~ (\neg\alpha_{p,s,m}(x)\rightarrow \neg \alpha_{p,s,m+1}(x))$ \quad  for all $p\vert n$, $s<k_p $ and $m \geq 1$
                \item[(AS4)] \label{Axiom:Sp4} $(\forall x)~ (\neg F(x) \rightarrow \neg \alpha_{p,s,1}(x))$\quad for all $p\vert n$ and $s<k_p$
                \item[(AS5)] \label{Axiom:Sp5} $(\forall x)~ (\neg A(x)\rightarrow \neg \beta_{p,1}(x))$\quad for all $p|n$
                \item[(AS6)] \label{Axiom:Sp6} $(\forall x)~ (A(x) \vee F(x))$\footnote{We adapted this axiom to our convention. Notice that in \cite{Sch82}, the corresponding axiom reads  $(\forall x)~ ((A(x) \vee F(x)) \wedge (A(x)\rightarrow \exists y (x<y))$. This seems to be a mistake: by their convention, the $n$-spine can contain a last element which is an $A$ (and an $F$). 
                }
                \item[(AS7)] \label{Axiom:Sp7}  $(\forall x)~  F(x) \leftrightarrow \bigvee\{\alpha_{p,s,1}(x) : p|n, s<k_p\})$\footnote{We adapted this axiom to our convention. In \cite{Sch82}, Schmitt needs to make an exception when there is a smallest convex subgroup which is divisible. }
                \item[(AS8)] \label{Axiom:Sp8} $(\forall x)~ (A(x)\wedge F(x)) \rightarrow[\bigwedge\{ \alpha_{p,k_{p}-1,m}(x) \leftrightarrow \beta_{p,m}(x) : p|n\}\wedge$
                \\
                $\bigwedge \{ \neg \alpha_{p,s,1}(x): p|n, ~ s<k_p \}  ]$\footnote{In \cite{Sch82},  the corresponding axiom has a small typo, as the $\neg$ is missing.}
                \item[(AS9)] \label{Axiom:Sp9} $(\forall x)~ (D(x) \rightarrow \beta_{p,1}(x) \wedge \neg \beta_{p,2}(x) ),$ for all $p|n$
                \item[(AS10)] \label{Axiom:Sp10} $(\forall x)~ \left[ A(x) \rightarrow  (F(x) \leftrightarrow \bigvee_{p\vert n} \beta_{p,1}(x))\right] $
                \item[(AS11)] \label{Axiom:Sp11} $(\forall x,y)~ \left[(x<y \wedge A(x)) \rightarrow (\exists z) ( x \leq z < y \wedge F(z) )\right] $
                \item[(AS12)]\label{Axiom:Sp12} $(\forall x,y)~ \left[(x<y \wedge \neg A(x)) \rightarrow (\exists z) ( x < z < y \wedge A(z) )\right] $
                \item[(AS13)] \label{Axiom:Sp13} $(\forall y) \left[\neg A(y) \rightarrow (\exists z)(z<y \wedge A(z))\right]$
                \item[(AS14)] \label{Axiom:Sp14} $(\forall x,y)~\Big[ (x<y \wedge F(y) \wedge \alpha_{p,s,1}(y)) \rightarrow (\exists z)( x \leq z < y \wedge F(z) \wedge \bigvee_{l<k_p} \alpha_{p,l,1}(z) \Big] $ 
                for all $p|n$ and $s< k_p-1$
            \end{enumerate}
        \end{definition}

        Here is one of the main result of Schmitt's \emph{Habilitationsschrift} (\cite[Lemma 5.4 \& Theorem 5.7]{Sch82}):
        \begin{fact}[Schmitt] \label{fact:Schmitt}
            For any ordered abelian group $G$, $\Sp_n(G)$ is a model of $T_{\Sp,n}$. Conversely, for all models $S$ of $T_{\Sp,n}$, there is an ordered abelian group $G$ such that $\Sp_ n(G) \equiv S$. 
        \end{fact}

        \begin{remark}\label{rem:endsegment}
        Let $M$ be a model of $T_{\Sp_n}$. Then an end-segment $M'$  automatically satisfies all axioms, except potentially the Axiom (AS13) (which asserts the existence of a point $y$ before an $\neg A$).   However, if $M'$ has no minimum, or if its minimum satisfies $A$, then $M'$ satisfies (AS13), and will therefore also be a model of $T_{\Sp_n}.$ 
        
        \end{remark}

        \begin{example}
    Consider $H\coloneqq \bigoplus_{\omega} \Zbb$ and let $a= (a_n)_{n\in \omega} \in \hprod_{\omega} \Zbb$ with $a_n=2$ for all $n\in \mathbb{N}$. Set 
    \[G \coloneqq \braket{H,a}.\]
    The $A_n$'s are exactly the following convex subgroups:
    \[G_k \coloneqq \{x\in G : \forall i<k ~x(i)=0\}.\]
    and we have $D(G_k) \wedge F(G_k)$ for all $k$ and all $n$. There is one additional $F_2$, and it is $F_2(a)=\{0\}$.
    
    \begin{center}
        \begin{tikzpicture}[scale=0.6]
            \foreach \j in {0,...,3}{
                    \filldraw ( \j,0) circle (1mm);
                    \draw (\j,0) node[anchor=north]{$G_\j$};
                };
            
            \filldraw[red] (6,0) circle (1.3mm);
            \draw (6,0) node[anchor=north]{$\{0\}$};
            \draw (4.5,0) node{$\cdots$};
            
        %}
        \end{tikzpicture}
    \end{center}
    The $2$-spine $\{G_n : n\in \mathbb{N} \} \cup \{ ~\{0\}~\}$, equipped with the reversion inclusion order and the appropriate colors $\alpha_{2,k,m},\beta_{2,m},D,F_2,A_2$, is a model of $T_{\Sp_2}$. Also, all end-segments are model of $T_{\Sp_2}$, except the end-segment consisting only of the last point.
    
\end{example}
 We give later a more involved example of a group $G$ which $2$-spine admits an augment which is not a model of (see \Cref{ex:DelonLucas}).

    \begin{lemma}[{\cite[Lemma 3.15]{Sch82}}]\label{lem:SpmtranslatetoSpn}
        Let $m,n$ be integers such that $m|n$. For every sentence $\phi$ in $\mathcal{L}_{\Sp_m}$, there is a sentence $\phi_n$ in $\mathcal{L}_{\Sp_n}$ such that, for all ordered abelian groups $G$, 
        \[\Sp_m(G) \models \phi \Leftrightarrow \Sp_n(G) \models \phi_n.\]
    \end{lemma}
    In fact, \cite[Lemma 3.15]{Sch82} is stronger and can be reformulated in terms of interpretation. The statement above is however all we need for the main result of this section.

    Unfortunately, we do not have, in general, an inclusion $\Sp_n \subseteq \Sp_m$ for $n\vert m$. In other word, the set of spines is not (quite) a directed system of chains:  
        \begin{example}

            Consider $G \coloneq \bigoplus_{\omega}\Ztwo \oplus \bigoplus_{\zeta}\Zthree$. We compute the order type of the $2$,$3$ and $6$-spines.
           Consider the elements $a=(1,0,0,\dots, \dots, 0,0,0, \dots)$ and $b=(0,0,0,\dots, \dots, 0,1,0, \dots)$. Then we have:
    \begin{itemize}
        \item $A_3(a)= \bigoplus_ {\zeta}\Zthree$ and $A_6(a)= \bigoplus_{\omega\setminus\{0\}}\Ztwo \oplus \bigoplus_ {\zeta}\Zthree$;
        \item $A_2(b)= \{0\}$ and $A_6(b)= \bigoplus_ {n>0}\Zthree$.
    \end{itemize}
     Then, one can see that $\Sp_2$ is of type $\omega +1 $ and contains  an extra point for $\{0\}$. The $6$-spine $\Sp_6$ is as expected  of type $\omega + \zeta$ but the spine $\Sp_3$ contains one extra point for the convex subgroup $ \bigoplus_ {\zeta}\Zthree$.
\begin{center}
    \begin{tikzpicture}[scale={0.6}]
    \foreach \i in {0,1,2}{
        \filldraw (\i,0) circle (1mm);
        \filldraw (\i,-4) circle (1mm);
    };
    \foreach \i in {6,7,8}{
        \filldraw (\i,-2) circle (1mm);
        \filldraw (\i,-4) circle (1mm);
    };
    \filldraw (4,-2) circle (1mm);
    \filldraw (10,0) circle (1mm);
  
    \draw (3,0) node {$\dots$};
    \draw (3,-4) node {$\dots$};
    \draw (5,-2) node {$\dots$};
    \draw (5,-4) node {$\dots$};
    \draw (9,-2) node {$\dots$};
    \draw (9,-4) node {$\dots$};
    %\filldraw (4,0) circle (1mm);   
    \draw[dashed] (3.6,-5) -- ++ (0,5.5);
    \draw[dashed] (9.6,-5) -- ++ (0,5.5);
    \draw[dashed] (10.4,-5) -- ++ (0,5.5);
    \draw[dashed] (4.4,-
    5) -- ++ (0,5.5);

    \draw (-2,0) node {$\Sp_2(G)$};
    \draw (-2,-2) node {$\Sp_3(G)$};
    \draw (-2,-4) node {$\Sp_6(G)$};
    \draw (4,-5) node[anchor=north] {$\bigoplus_\zeta \Zthree$};
    \draw (10,-5) node[anchor=north] {$\{0\}$};
    \end{tikzpicture}
\end{center}
    In other words, $2$-divisible ribs can ``create'' a point in the $3$-spine and \emph{vice versa}.
\end{example}
%    \begin{lemma}\label{lem:SpnInterprSpm}
 %       Let $m,n$ be integer such that $m|n$. There is an interpretation $\iota$ of $T_{\Sp,m}$ in $T_{\Sp,n}$ such that for all ordered abelian groups,
 %       $\iota (\Sp_n(G)) \cong \Sp_m(G)$.
 %   \end{lemma}

\subsection{Weak Augmentability by Infinitesimals}

    The next example is due to Delon and Lucas. Rephrased in our setting, it is an ordered abelian group whose $2$-spine $\Sp_2(G) \models T_{\Sp_2}$ can be ``augmented'', i.e. there is a $\mathcal{L}_{\Sp_2}$-structure $\Sigma$ such that  $\Sp_2(G) + \Sigma \models T_{\Sp_2}$, but the augment $\Sigma$ itself doesn't satisfy (AS13).

        \begin{example}[{\cite[Example p.503-4]{DelonLucas}}]\label{ex:DelonLucas}
    Consider $H\coloneqq \bigoplus_{\omega^\omega} \Zbb$ and $H'\coloneqq \bigoplus_{\omega^\omega\cdot 2} \Zbb$. Let $(\alpha_i)_{i<\omega}$ be a sequence of ordinals cofinal in $\alpha$.  Let $a^\alpha= (a^\alpha_\beta)_{\beta\in \omega^\omega} \in \hprod_{\omega^\omega} \Zbb: $
    \[a^\alpha_\beta=\begin{cases}
         2 & \text{ if } \beta=\alpha_i \text{ for some }i<\omega,\\
        0 & \text{ otherwise.}
    \end{cases}\]

    Then set 
    \[G \coloneqq \braket{H,a^\alpha: \alpha<\omega^\omega, \text{ limit}}\]
    and
    \[G' \coloneqq \braket{H',a^\alpha: \alpha<\omega^\omega\cdot 2, \text{ limit}}.\]
    The convex subgroups are the
    \[G_\alpha \coloneqq \{x\in G : \forall \beta>\alpha ~x(\beta)=0\}\]
    where $\alpha \in \omega^\omega\cdot 2$.
    The convex subgroups $G_{\alpha+1}$ (successor ordinals) are $D$ and $A_n,F_n$ for all $n$. All $G_{\alpha}$ (successors and limits) except $G_0 = G$ are $F_2$'s. We have also for $\alpha$ limits $\alpha_{2,s}(\Gamma(n,G_{\alpha})) = 
    \begin{cases}
        1 & \text{ if } s=1 \text{ and } 2|n\\
        0 & \text{ otherwise. }
    \end{cases}$
    \begin{center}
        \begin{tikzpicture}[scale=0.4]
            \foreach \i in {0,...,3}{
                \foreach \j in {1,...,3}{
                    \filldraw ( 6*\i+\j,0) circle (1mm);
                };
                \filldraw[red] (6*\i,0) circle (1.3mm);
                \draw (6*\i+4,0) node{$\cdots$};
            }
            \filldraw[white] (0,0) circle (1.5mm);
            \draw (6,0) node[anchor=north]{$G_\omega$};
            \draw (12,0) node[anchor=north]{$G_{\omega\cdot2}$};
            \draw (18,0) node[anchor=north]{$G_{\omega\cdot3}$};
            
        \end{tikzpicture}
    \end{center}
    We consider the end-segment $\Sigma \coloneqq\{G_\alpha :  \omega^\omega \leq \alpha <\omega^\omega \cdot 2 \}$
    \begin{center}
        \begin{tikzpicture}[scale=0.4]
            \foreach \i in {0,...,3}{
                \foreach \j in {1,...,3}{
                    \filldraw ( 6*\i+\j,0) circle (1mm);
                };
                \filldraw[red] (6*\i,0) circle (1.3mm);
                \draw (6*\i+4,0) node{$\cdots$};
            }

            \draw (0,0) node[anchor=north]{$G_{\omega^\omega}$};
            \draw (6,0) node[anchor=north]{$G_{\omega^\omega+\omega}$};
            \draw (12,0) node[anchor=north]{$G_{\omega^\omega+\omega\cdot 2}$};
            \draw (18,0) node[anchor=north]{$G_{\omega^\omega+\omega\cdot 3}$};
            
        \end{tikzpicture}
    \end{center}
    We can see that 
    \[\Sp_2{G} + \Sigma = \Sp_2{G'}  \equiv \Sp_2{G},\]
    but $\Sigma$ is not the $2$-spine of a group, as its doesn't satisfy (AS13).
\end{example}

\begin{fact}[Delon-Lucas {(\cite[Proposition 4-4 2.]{DelonLucas})}]\label{fact:SpineSumOAGS}
    Let $G$ and $H$ be two ordered abelian groups and let $n$ be an integer. Assume that one of these two statements \emph{does not} hold:\footnote{Note that Delon-Lucas use the inclusion-preserving-order on the spines, which must be reversed to fit our convention.}
    \begin{enumerate}
        \item $\Sp_n(G)$ has a last element $x$ which satisfies $\neg F_n$.
        \item $\Sp_n(H)$ has a first element $y$.
    \end{enumerate}
    Then $\Sp_n(G \oplus H)=\Sp_n(G) + \Sp_n(H)$.
\end{fact}

    \begin{definition}
        Let $G$ be an ordered abelian group. We denote $\Sp(G)$ the set $\{ \Sp_n(G) : n \geq 2\}$ of spines of $G$.
        We see $\Sp(G)$ as a multisorted structure, with a sort for each $\Sp_n(G)$ in the language $\lL_n$ and no relation nor function between them. 
    \end{definition}
    \begin{definition}

        We say that $\Sp(G)$ is weakly augmentable on the right if there is an integer $k\geqslant 1$ and a set of non-empty coloured chains $\Sigma \coloneqq \{ \Sigma_n : k\vert n\}$, $\Sigma_n$ an $\lL_n$-structure, such that
        \[ \Sp_n(G)+ \Sigma_n  \equiv \Sp_n(G)
        .\]    \end{definition}

        \begin{lemma}\label{lem:CharNonAugmByQ}
        Let $G$ be an ordered abelian group. TFAE:
        \begin{enumerate}
        \item $G$ is not weakly augmentable by $\mathbb Q$.
        \item there exists $k$ such that $\Sp_k(G)$ either has no final point or has a final point satisfying $F_k$.
        \item there exists $k$ such that for all $n$ multiple of $k$, $\Sp_n(G)$ either has no final point or has a final point satisfying $F_n$.
        \end{enumerate}
        \end{lemma}

\begin{proof}
 $(1)\Leftrightarrow(2)$ comes from \Cref{DOAG-aug}-\p0. To prove $(2)\Leftrightarrow(3)$, we simply observe that for any $n$, $\Sp_n(G)$ has a final point satisfying $\neg F_n$ if and only if $G$ has a non-trivial $n$-divisible convex subgroup. This is in fact part of the proof of \Cref{DOAG-aug}.
 
 %Indeed, if $G$ has a non-trivial $n$-divisible convex subgroup, name it $C$, then any element $c\in C$ is such that $A_n(c)=\{0\}$, so $\Sp_n(G)$ has a final point. For any $g\in G$, either $g\in nG$ and $F_n(g)=\emptyset$ is not in the spine, or $g\notin nG$, but then $g$ is not $n$-divisible modulo $C$, so $F_n(g)\supseteq C$ and is not the final point of $\Sp_n(G)$. Conversely, if $\Sp_n(G)$ has a final point satisfying $\neg F$, then this point satisfy $A$, so there is $g\in G$ such that $A_n(g)=\{0\}$. Now $B_n(g)$ is a convex subgroup of $G$. Take $h\in B_n(g)$, then $F_n(h)\subseteq A_n(h)=A_n(g)=\{0\}$ by properties of the spine, since the final point is not an $F$, we must have $F_n(h)=\emptyset$, that is, $h$ is $n$-divisible.

Now $(3)\Rightarrow(2)$ is immediate, and we also have the implication $(2)\Rightarrow(3)$, since if no convex subgroup of $G$ is  $k$-divisible, then no convex subgroup of $G$ is $mk$ divisible, for any $m\in\mathbb N$.
\end{proof}
    
    \begin{theorem}\label{thm:WeakAugmentableInfinitesimalsOags}
        Let $G$ be an ordered abelian groups. If $\Sp(G)$ is weakly augmentable on the right, then $G$ is weakly augmentable by infinitesimals.
        If $G$ is not weakly augmentable by $\Qbb$, then the converse also holds: if $G$ is weakly augmentable by infinitesimals, then 
        $\Sp(G)$ is weakly augmentable on the right.
    \end{theorem}

    \begin{proof}
        To show the first implication, assume that $\Sp(G)$ is weakly augmentable on the right: let $k\in \mathbb{N}$ and  $\{\Sigma_n: k\vert n\}$ all non-empty be such that
        \[ \Sp_{n}(G)+ \Sigma_n  \equiv \Sp_{n}(G)\]
        for all integers $n$ such that $k\vert n$.

        If $\Sigma_n$ has a first element $y$ satisfying $\neg A_n(y)$, then from (AS12) applied to $\Sp_{n}(G)+ \Sigma_n$ at the point $y$, we deduce that $\Sp_{n}(G)$ has no maximum element. In such a case, $\Sp_{n}(G)$ is in fact strongly augmentable (\Cref{lem:aug-mo}). Then, by \Cref{lem:StrongAugment}, we may find a strong augment of $\Sp_{n}(G)$ with no minimal element.
        In any case, we may assume that $\Sigma_n$ doesn't have a first element satisfying $\neg A_n(y)$.
        
        By \Cref{rem:endsegment}, $\Sigma_n$ satisfies $T_{\Sp_n}$ and by \Cref{fact:Schmitt}, there are groups $H_n$ with $\Sp_n(H_n)\equiv \Sigma_n$.

        We show now that $G$ and $H_n$ satisfy the conditions of \Cref{fact:SpineSumOAGS} (i.e. either $(1)$ or $(2)$ is false), and thus, we would deduce that:
        \[\Sp_n(G\oplus H_n) =\Sp_n(G) + \Sp_n(H_n) = \Sp_n(G) + \Sigma_n \equiv \Sp_n(G). \]

        Indeed, the problematic case is when $\Sp_n(G)$ has a last element $x$ satisfying $\neg F_n(x)$ and $\Sigma_n$ has a first element $y$. This is impossible, as such an $x$ would then satisfy $A_n(x)$ by (AS6). But now, by (AS11), as $y$ is the immediate successor of $x$, $x$ would need to satisfy $F_n$, a contradiction.

        Let $H \coloneq \prod_\Ucal H_n$ be a non-principal ultraproduct of the $H_n$, with $n\mathbb{N} \in \Ucal $ for all $n$. Let $m$ be an integer with $k\vert m$ and let $\phi$ be a sentence in the language of $m$-spines. For all multiple $n$ of $m$, by \Cref{lem:SpmtranslatetoSpn}, there is a sentence $\phi_n$ in the $n$-spine which interprets $\phi$. Also, observe that $G \oplus H \equiv (\prod_\mathcal{U}G) \oplus (\prod_\mathcal{U}H_n) \equiv \prod_\mathcal{U}(G \oplus H_n)$, and therefore, by Łoś's theorem:
        \begin{align*}
            \Sp_m(G \oplus H) \models \phi &\Leftrightarrow  \Sp_m(\prod_\mathcal{U}(G \oplus H_n)) \models \phi,\\
            &\Leftrightarrow \text{for $\mathcal{U}$-many $n$,} \Sp_m(G \oplus H_n) \models \phi,\\
            &\Leftrightarrow \text{for $\mathcal{U}$-many $n$ with $m|n$,} \Sp_n(G \oplus H_n)= \Sp_n(G)+\Sigma_n \models \phi_n \\
            &\Leftrightarrow  \text{for $\mathcal{U}$-many $n$ with $m|n$,} \Sp_n(G) \models \phi_n\\
            &\Leftrightarrow  \text{for $\mathcal{U}$-many $n$ with $m|n$,} \Sp_m(G) \models \phi\\
            &\Leftrightarrow \Sp_m(G) \models \phi
        \end{align*}
            
        By \Cref{prop:SpinesInterpretableInOAG}, $G \oplus H \equiv G$, as wanted.
        
       To show the second implication, assume that there is $H$ such that $G \oplus H\equiv G$, and $G$ is not augmentable by some divisible ordered abelian group. By \Cref{lem:CharNonAugmByQ}, there is an integer $k$ such that, for all multiple $n$ of $k$, $\Sp_n(G)$ doesn't admits a last point satisfying $\neg F_n$. By \Cref{fact:SpineSumOAGS},
       \[\Sp_n(G \oplus H) =\Sp_n(G)+ \Sp_n(H).\]
       By \Cref{prop:SpinesInterpretableInOAG}, $\Sp_n(G) \equiv \Sp_n(G) + \Sp_n(H) $ for all such $n$. Since $H$ is not divisible, $\Sp_n(H)$ is not empty, in particular, $\Sp(G)$ is augmentable on the right, as wanted.
    \end{proof}

    Merging this result with \Cref{DOAG-aug}, we obtain:

    \begin{corollary}
    An ordered abelian group is weakly augmentable by infinitesimals if and only if one of the following holds:
    \begin{enumerate}
    \item it is weakly augmentable by a divisible infinitesimal augment, in which case it is also strongly augmentable by an infinitesimal augment,
    \item its multisorted set of spines is augmentable.
    \end{enumerate}
    \end{corollary}

    Note that the two cases are not mutually exclusive, as seen in \Cref{ex:aug-by-Q-and-Z+Q}.

    We conclude this section with simple examples, to show that for all integers $n\neq m$, there is no correlation between right-augmentability of $\Sp_n(G)$ and right-augmentability of $\Sp_m(G)$.
    \begin{example}
        Denote by $\eta$ the linear order $(\mathbb{Q},<)$. Consider the ordered abelian group 
        \[G \coloneq \bigoplus_{\eta}(\Ztwo \oplus \Qbb) \oplus (\Ztwo \oplus \Zthree).\] 
        We have
        \[\Sp_2(G)=(\eta+1) \cdot L\]
        where $L= \{a,b\}$ with $A_2(a)\wedge F_2(a)$ and $A_2(b)\wedge \neg F_2(b)$. In particular, $\Sp_2(G)$ is right-augmentable and in fact $\Sp_2(G)+\Sp_2(G) \cong \Sp_2(G)$.
        The $6$-spine $\Sp_6(G)$ is $\eta \cdot L + L'$ where 
        \begin{itemize}
            \item $L= \{a,b\}$ with $A_6(a)\wedge F_6(a)$ and $A_6(b)\wedge \neg F_6(b)$.
            \item $L'= \{a,b'\}$ with $A_6(a)\wedge F_6(a)$ and $A_6(b')\wedge  F_6(b')$. 
        \end{itemize}
        We can see that the $6$-spine is not right-augmentable and in fact, $G$ is not augmentable by infinitesimals. The converse is technically also false:
    \end{example}
    \begin{example}
        Consider $G\coloneqq \bigoplus_{\omega*}  \Ztwo$. The $6$-spine is right-augmentable, but the $3$-spine is a singleton-chain $\{a\}$ (with $\neg F_3(a)\wedge A_3(a)$ by convention) and therefore is not augmentable.
    \end{example}

        In the next example, we give an ordered abelian group $G$ such that $\Sp_{p^k}(G)$ is augmentable on the right, but $\Sp_{p^{k+1}}(G)$ is not. This is due to our convention to not include all predicates $\alpha$ in the spines.\footnote{For example, according to Schmitt's definitions that we use here, only the predicates $\alpha_{p,1,m}$, $m\geq 1$, are present in the $p$-spine. In \cite{CH11} however, the $p$-spine is equipped with all the predicates $\alpha_{p,s,m}$, $s\geq 1,m \geq 1$.}
    \begin{example}
        Consider the ordered abelian groups $H= \braket{\bigoplus_{\omega} \mathbb{Z}, a}$ and $H'= \braket{\bigoplus_{\omega} \mathbb{Z}, b}$, where $a= (2)_{n\in \omega}$ and $b=(2^2)_{n\in \omega}$ are elements of $\hprod_{\omega}\mathbb{Z}$.
        Set $G= \bigoplus_{\omega^*} H \oplus H'$. 
        Then $\Sp_4(G)$ is augmentable on the right but $\Sp_8(G)$ is not.
        
        The $4$-spine of $G$ can indeed be identified with the linear order $\bigoplus_{\omega^*}(\omega+1)$. The convex subgroups $C$ corresponding to the points $(n,m)$, $n,m\in \omega$ have colors $D(C)$ and the convex subgroups $C$ corresponding to the points $(n,\omega)$, $n\in \omega$ have colors $F_2(C)\wedge \neg A_2(C) \wedge \alpha_{2,1,1}(C)\wedge \neg \alpha_{2,1,2}(C)$.
   \begin{center}
        \begin{tikzpicture}[scale=0.4]
            \foreach \i in {3,...,0}{
                \foreach \j in {1,...,3}{
                    \filldraw ( 6*\i+\j,0) circle (1mm);
                };
                \filldraw[red] (6*\i+5,0) circle (1.3mm);
                \draw (6*\i+4,0) node{$\cdots$};
            }
        \draw (-2,0) node{$\cdots$};
        \draw (17,0) node[anchor=north]{$(1,\omega)$};
        \draw (11,0) node[anchor=north]{$(2,\omega)$};
        \draw (23,0) node[anchor=north]{$(0,\omega)$};
            
        \end{tikzpicture}
    \end{center}

        The $8$-spine of $G$ is identical, except the last point $(0,\omega)$ has colors $F_2(C)\wedge \neg A_2(C)\wedge \alpha_{2,2,1}(C)\wedge \neg \alpha_{2,2,2}(C) \wedge \neg \alpha_{2,1,1}(C)$.

       \begin{center}
        \begin{tikzpicture}[scale=0.4]
            \foreach \i in {3,...,0}{
                \foreach \j in {1,...,3}{
                    \filldraw ( 6*\i+\j,0) circle (1mm);
                };
                \filldraw[red] (6*\i+5,0) circle (1.3mm);
                \draw (6*\i+4,0) node{$\cdots$};
            }
        \draw (-2,0) node{$\cdots$};
        \draw (17,0) node[anchor=north]{$(1,\omega)$};
        \draw (11,0) node[anchor=north]{$(2,\omega)$};
        \filldraw[black!40!green] (23,0) circle (1.6mm);
        \draw (23,0) node[anchor=north]{$(0,\omega)$};
            
        \end{tikzpicture}
    \end{center}
    \end{example}

%auto-ignore
\section{Automatic Definability of henselian valuations}\label{sec:AutomaticDefinability}

We are now ready to characterize automatic definability of valuations. Recall that an ordered abelian group is said to have automatic ($\emptyset$)-definability (in the class of henselian valued fields) if for any henselian valued field $(K,v)$ with $vK\cong G$, the valuation $v$ is $\Lring$-($\emptyset$)-definable, and similarly for a (residue) field $k$.

\subsection{By properties of the residue field}

We build on \Cref{thm:strongleftaugmentabilityOAG} to study elementary embeddings of henselian valued fields, and to characterize when the valuation is automatically (uniformly) definable from properties of the residue field.

Recall that a field $k$ is called $t$-henselian if $k$ is $\mathcal{L}_\textrm{ring}$-elementarily equivalent to a field admitting a non-trivial henselian valuation.
The key observation in what follows is that \Cref{thm:strongleftaugmentabilityOAG} gives rise to an elementary embedding of any $t$-henselian field of characteristic $0$ into some (non-trivial) power series field over it:

\begin{proposition}
    Let $k$ be a $t$-henselian field of characteristic $0$. Then 
    $k \substruct k((H))$ for some non-trivial ordered abelian group $H$. \label{power}
\end{proposition}
\begin{proof}
    Let $k$ be a $t$-henselian field of characteristic $0$. 
    Arguing as in \Cref{rem:soft}, we may replace $k$ with any $\mathcal{L}_\textrm{ring}$-elementary extension.
In particular, 
    we may assume that $k$ is $\aleph_1$-saturated in $\Lring$.
    %$k$ admits a nontrivial henselian valuation.
    By \cite[Lemma 3.3]{PZ78}, $k$ admits a non-trivial henselian valuation $v$. 
    Passing to an $\aleph_1$-saturated $\mathcal{L}_\textrm{val}$-elementary extension $(k^*,v^*)$ of $(k,v)$, we may further assume that $(k,v)$ is $\aleph_1$-saturated as a valued field.
    We next argue that $k$ admits a non-trivial henselian valuation $w$
    of residue characteristic $0$. 
    Note that $v$ has a non-trivial coarsening $w$ of residue characteristic $0$ if and only if 
    $vk\supsetneq \CVX \langle v(p) \rangle$ holds. That this inclusion is proper is ensured by $\aleph_1$-saturation of $(k,v)$. Thus, replacing $v$ by $w$, we may assume $\mathrm{char}(kv)=0$.

    Now, by \Cref{fact:EmbedsHahnSeries}, we have $(k,v)\substruct (kv((vk)),u)$, where $u$ denotes the power series
    valuation on $kv((vk))$. %     This argument is also in \cite{vdD14}, and we include     it for the convenience of the reader.      Since  $(k,v)$ be a $\aleph_1$-saturated, there is a section of the valuation $s:vk \rightarrow k$. Since $(k,v)$ is of equicharacteristic $0$, there is a lift $l:kv\rightarrow k$ of the residue field. In particular, it follows that $kv(vk)$ embeds in $k$. The extension is immediate and therefore, $k$ and $kv(vk)$ admit $kv((vk))$ for unique (up to isomorphism) maximal extension. Then we have  $(k,v) \substruct (kv((vk)),u)$ by Ax-Kochen/Ershov.
Denote $G:=vk$. By Theorem \ref{thm:strongleftaugmentabilityOAG}, we have
    $G \substruct H \oplus G$ for some non-trivial ordered abelian group $H$.
    Thus, we obtain -- once again by the Ax-Kochen/Ershov Theorem in equicharacteristic
    $0$ --
    $$k \substruct kv((G)) \substruct kv((G))((H))$$
    and $$k((H)) \substruct kv((G))((H)).$$
    This induces an elementary embedding $k \substruct k((H))$.
\end{proof}

We now link the proposition above to $\Lring$-definability of henselian
valuations with a given residue field. 

\begin{theorem} \label{def_fieldwise}
    Let $k$ be a field of characteristic $0$. The following are equivalent:
    \begin{enumerate}
        \item $k$ is not $t$-henselian, i.e., it is not elementary equivalent to a field $k'$ that admits a non-trivial henselian valuation,
        \item $k$ has automatic definability, i.e., for every henselian valued field $(K,v)$ with residue field $k$, the valuation $v$ is definable (possibly using parameters),
        \item $k$ has automatic $\emptyset$-definability, i.e., for every henselian valued field 
        $(K,v)$ with residue field $k$, the valuation $v$ is 
        $\emptyset$-definable,
        \item all henselian valuations with residue field elementarily equivalent to $k$ are
        uniformly $\emptyset$-definable in $\Lring$. More precisely,
        there is a parameter-free $\Lring$-formula $\psi(x)$ which defines the valuation ring $v$ of in any henselian valued field $(K,v)$ with residue field $Kv$ a model of $\mathrm{Th}_{\Lring}(k)$.
    \end{enumerate}
\end{theorem}
\begin{proof} The implications (4) $\implies$ (3) and (3) $\implies$ (2) are trivial.

    We show (2) $\implies$ (1) via contraposition.
    Assume that $k$ is $t$-henselian. By Proposition \ref{power},
    we have $k \substruct k((H))$ for 
    some non-trivial ordered abelian group $H$. 
    We have the following chain of elementary embeddings:
    \[k \substruct k((H_0)) \substruct   k((H_{-1} \oplus H_0 \oplus H_{1})) \substruct \cdots \substruct k((\bigoplus_{-n<i<n} H_i)) \substruct \cdots\]
    where $H_i$ are all copies of $H$ and with the obvious embeddings. By Tarski's chain lemma \cite[Theorem~2.1.4]{TZ12}, $k$ is an elementary substructure of $k':=k((\bigoplus_{\zeta}H))$. By \Cref{thm:straugvf}, it is enough to show that $k'$ is strongly co-augmentable with some ordered abelian group $\Gamma$. %show that there is a henselian valued field with residue $k'$ whose valuation is not definable in $\Lring$. The same property will hold for $k$ instead of $k'$ by the $\substruct$-version of the      Ax-Kochen/Ershov theorem (\cite[Theorem 6.17]{Hils-mtvf}).
    Set  $\Gamma:=\bigoplus_{\zeta} H$ and let %consider the Hahn field    $K=k'((\Gamma))$, together with the power series    valuation $v$ with value group $\Gamma$. Consider the field    $K^* = k'((\Gamma'))((\Gamma))$ with
    $\Gamma'$ be another copy of $\Gamma$.  
   Note that $\Gamma \substruct \Gamma\oplus\Gamma'$ 
   (by, e.g., playing Ehrenfeucht–Fraïssé games) and that $k' \substruct k'((\Gamma'))$. In other words, $k'$ and $\Gamma$ are strongly co-augmentable by $\Gamma'$, as claimed.
   %The field $K^*$ admits the distinct henselian valuations   \begin{itemize}  \item $u$ with residue field $k'$ and value group $\Gamma \oplus \Gamma'$ and        \item $w$ with residue field $k'((\Gamma'))$ and value group $\Gamma$.    \end{itemize}
   %Then we have by applying Ax-Kochen/Ershov \cite[Theorem 6.17]{Hils-mtvf}
   %once again
   %$$(K,v)\substruct (K^*,u) \textrm{ and }(K,v)\substruct (K^*,w)$$
%   for the natural (field) embedding of $K$ in $K^*$.
  %  By Corollary \ref{comparable}, this shows that $v$ is not $\mathcal{L}_\textrm{ring}$-definable on $K$.

    The implication (1) $\implies$ (4) is essentially \cite[Proposition 5.5]{FJ15}. We give the proof for convenience of the reader.
    Assume $k$ is not $t$-henselian, in particular $k$ is not separably closed.
    Consider the $\mathcal{L}_\textrm{val}$-theory $T$ which stipulates that any model $(K,v)$ is henselian with residue field $Kv \equiv k$ (in $\mathcal{L}_\textrm{ring}$) and the 
    $\mathcal{L}_\textrm{val}$-formula $\phi(x)$ which asserts that $ x \in \mathcal{O}_v$.
    Take $(K_1,v_1), (K_2,v_2) \models T$ that have the same $\Lring$-reduct, i.e., in particular $K_1=K_2$. Then both $v_1$ and $v_2$ are henselian valuations with
    non-separably closed residue field, thus they are comparable, cf.~\cite[Theorem 4.4.2]{EP05}. 
    Assume for a contradiction that we have $\mathcal{O}_{v_1} \subsetneq \mathcal{O}_{v_2}$, then $v_1$ induces a
    non-trivial henselian
    valuation on $Kv_2$. But then $k$ is $t$-henselian, in contradiction to our assumption. Thus,
    we have $v_1=v_2$, and so \Cref{fact:beth} implies that $v$ is uniformly $\emptyset$-definable in
    all models of $T$.
\end{proof}

\begin{remark} \Cref{def_fieldwise} does not hold in general for fields of positive characteristic.
\begin{enumerate}
\item It is always true, regardless of the characteristic, that $k$ not being $t$-henselian implies $k$ has 
automatic $\emptyset$-definability, see \cite[Proposition 5.5]{FJ15}. In fact, the argument in \cite{FJ15} even shows that the valuation ring is $\emtpyset$-definable by an $\exists\forall$-formula.
\item The converse direction fails in general. If $k$ is imperfect and admits no Galois extensions of degree divisible by $p$, then any henselian valuation with residue field $k$ is definable by \cite[Proposition 3.6]{Jah24}. Such a field could be $t$-henselian and still have automatic definability. Examples exist abundantly, e.g., $\mathbb{F}_p(t)^\mathrm{sep}$ is $t$-henselian (and in fact even henselian) but admits no proper
Galois extensions of any degree.
\item We do not know whether \Cref{def_fieldwise} hold for perfect fields of positive characteristic, and we have no counterexample. A crucial point in our proof is \Cref{power}, which is false for perfect fields in general: if $k$ is perfect, then $k \substruct k((H))$ would imply that $H$ is $p$-divisible and hence the power series valuation on $k((H))$ is tame.
%By [Jonas MSc thesis], there are perfect henselian fields $k$ which
However, there are examples of perfect henselian fields which have no elementary extension
that admits a non-trivial tame
valuation.
One such example is
$k=(\mathbb{F}_p((t))^\mathcal{U})^\mathrm{perf}$, the perfect hull of any nonprincipal
ultrapower of the power series field $\mathbb{F}_p((t))$ (see
\cite[Example 4.1.5]{JonasMSc}).
\item It is nonetheless possible to prove \Cref{power} for some classes of fields of positive characteristic; namely, if $k$ is $t$-henselian \emph{of tame type}, that is, if $k$ is elementary equivalent to some field $L$ admitting a non-trivial henselian tame valuation $v$, then $k\substruct k((H))$ holds for some non-trivial ordered abelian group $H$. The proof of this is identical to the proof of \Cref{power}, using the tame version of Ax-Kochen/Ershov (\cite[Theorem 7.1]{Kuh16}) and noting that since the value group $vL$ is $p$-divisible, any infinite augment $H$ of $vL$ must also be $p$-divisible.

Conversely, if $k$ is perfect and $k \preccurlyeq k((H))$ holds, then $k$ is necessarily $t$-henselian of tame type.
\end{enumerate}
\end{remark}

As we mentioned in the introduction, \Cref{def_fieldwise} partially answers the question raised in \cite[Section~5.4]{KKL}. Specifically, they inquire about identifying the classes $\Ccal$ and $\Ccal_0$ of fields with automatic definability and automatic $\emptyset$-definability respectively. Theorem \ref{def_fieldwise} shows that restricting to the case of characteristic $0$,
the class of fields that are not $t$-henselian coincides
exactly with both the classes $\Ccal$ and $\Ccal_0$. These classes are not elementary:

\begin{remark}\label{rem:not-axomatisable}
The class of $t$-henselian fields is axiomatisable in $\Lring$: it is clearly closed under elementary equivalence, and if $(K_i)_{i\in I}$ is a family of $t$-henselian fields, then for each $i$ there is a field $L_i\equiv K_i$ admitting a nontrivial henselian valuation $v_i$. Now, we have $\prod_{\mathcal U} K_i\equiv \prod_{\mathcal U} L_i$, and passing to $\Lval$, we have that $\prod_{\mathcal U} (L_i,v_i)$ is a henselian non-trivially valued field for any nonprincipal ultrafilter $\mathcal{U}$. Thus, $\prod_{\mathcal U} K_i$ is $t$-henselian, and the class of $t$-henselian fields is closed under ultraproducts.

However, its complement is not axiomatisable: calling a valued field
$n_\leq$-henselian if Hensel's Lemma holds for polynomials of degree up to $n$,
one can construct $n_\leqslant$-henselian fields
$(K_n,v_n)$ for every $n \geqslant 2$ which are not $t$-henselian
(see e.g.~\cite[Lemma 6.4]{FJ15}). Any non-principal ultraproduct of these is henselian.

Therefore, \Cref{def_fieldwise} implies that the classes $\Ccal$ and $\Ccal_0$ are not axiomatisable.
\end{remark}

\subsection{By properties of the value group}

We finally characterize value groups with automatic definability and automatic $\emptyset$-definability. The two notions are known to be distinct: in \cite[Example 4.9]{KKL} shows that $\bigoplus_{\omega^*}\mathbb Z$ has automatic definability but does not have automatic $\emptyset$-definability. 
Let us point out that this group is weakly augmentable by the infinitesimal $\mathbb Z$ but not strongly augmentable by infinitesimals, as seen in \Cref{ex:OAGNotStrongAug}.

% \begin{example}\pierre{I think this example has to be reworked a little. We have now \Cref{ex:OAGNotStrongAug}, we should make some reference}\blaise{I think we should just say: ``automatic def and auto 0 def are known to be different, see KKL ex. 4.9 for an example of a value group with auto def but not auto 0 def.''}
%  Take $G=\bigoplus_{\omega^*}\mathbb Z$, as in \Cref{ex:OAGNotStrongAug}. Since $G$ has a minimum positive element, one can define any henselian valuation with value group $G$ using Robinson's formula, with a parameter for an element with minimum positive valuation. That is, $G$ has automatic definability with parameters.
%  Now let $k=\mathbb C((\bigoplus_{\omega}\mathbb Z))$. $k$ and $G$ are weakly co-augmentable by $\mathbb Z$, so the valuation in $k((G))$ is not $\emptyset$-definable. That is, $G$ doesn't have automatic $\emptyset$-definability.
% \end{example}

\subsubsection{With parameters}

We give a complete characterization of ordered abelian groups with automatic definability, answering \cite[Question 5.1~(i)~($\dagger$)]{KKL}.

\begin{theorem}\label{str-abc}
 An ordered abelian group $G$ has automatic definability with parameters if and only if $G$ is not strongly augmentable by infinitesimals.
\end{theorem}

Recall that we have classified strong augmentability in \Cref{thm:StronglyAugmentableInfinitesimalsOAGs}. This classifcation is expressed by conditions on the $n$-fundaments of $G$, denoted by $F_n(g)$ for $g\in G$ (see \Cref{def:spines}).

In order to prove the theorem and define the valuation in any henselian valued field with non-augmentable value group, we will use ``Robinson's generalised formula''.

\begin{proposition}\label{grf}
Let $(K,v)$ be a henselian valued field. Suppose there is $t\in\mathcal O_v$ such that $v(t)$ is not $p$-divisible for some prime $p$. If $p\neq\Char(Kv)$, then consider Robinson's formula
$$\phi_p(x,t)\colon (\exists y)(1+tx^p=y^p).$$
%defines the set $I_t=\sset{a\in K}{pv(a)>-v(t)}$. 
If $p=\Char(Kv)$, consider instead the formula
$$\phi_p(x,t)\colon(\exists y)(1+tx^p=y^p-y).$$
In either case, the $\Lring$-formula
$$\psi_p(x,t)\colon (\forall y)(\phi_p(y,t)\rightarrow\phi_p(xy,t))$$
defines the valuation ring 
\[\mathcal O_w=\sset{a\in K}{v(a)>0\vee v(a)\in F_p(v(t))}\]
which is a nontrivial coarsening of $v$.
In particular, if $F_p(v(t))=\{0\}$, then $\psi_p(x,t)$ defines $\mathcal O_v$.
\end{proposition}

This is a quite well-known generalization of Robinson's formula, which we learned
from
Scanlon, see for example \cite[Prop.~3.6]{Jah24} or \cite[Prop.~1.4.1]{BoiPHD}.

\begin{proof}
We first claim that $\phi_p(x,t)$ defines the set $I_t\coloneq\sset{a\in K}{pv(a)>-v(t)}$.
This is a standard consequence of Hensel's lemma:
If  $1+tx^p$ is a $p^{th}$-power, then $v(1+tx^p)$ is clearly $p$-divisible. 
If  $1+tx^p$ is of the form $y^p+y$, then $v(y)<0$ (as otherwise $v(1+tx^p)= v(y)\geq 0=v(1)$ which implies $v(tx^p)=0$, contradicting that $v(t)$ is not $p$-divisible). Then again, $v(1+tx^p)$ is $p$-divisible. Since $v(tx^p)=pv(x)+v(t)$ is not $p$-divisible by choice of $t$, this is possible only if $pv(x)+v(t)=v(tx^p)>v(1)=0$, that is, $pv(x)<-v(t)$.

Conversely, assume that $pv(x)>-v(t)$. Then $v(1+tx^p)=0$. If ${\Char(Kv)\neq p}$, then $1\in Kv$ is a single-root of the $Kv$-polynomial $Y^p-\res(1+tx^p)$. By Hensel's lemma, there is $y\in K$ such that $y^p=1+tx^p$, as wanted.
Similarly, in residue characteristic $\Char(k)= p$,  $1\in Kv$ is a single-root the $Kv$-polynomial  $Y^p+Y-\res(1+tx^p)$ and Hensel's lemma gives $y\in K$ such that $y^p+y=1+tx^p$.\\

Now observe that the formula $\psi_p(x,t)$ defines the multiplicative stabilizer of $I_t$, that is, the set
\[S_t=\sset{b\in K}{bI_t\subseteq I_t}.\]
We aim to show that $b\in S_t$ if and only if $v(b)>0$ or $v(b)\in F_p(v(t))$.

If $v(b)\geqslant 0$, we take $a\in I_t$, and clearly $pv(ba)>pv(a)>-v(t)$, that is, $b\in S_t$. If $v(b)\in F_p(v(t))$ and $v(b)<0$, we take $a\in I_t$. Recall that, by \Cref{lem:FnSChmitt}(1), $v(b)\in F_p(v(t))$ means that $v(t)+h$ is not $p$-divisible for all $h\in vK$ with $|h|<p  |v(b)|$. Assume towards a contradiction that $ab\notin I_t$. Then we have $pv(a)>-v(t)$ but $p(v(a)+v(b))<-v(t)$. So, $0<pv(a)+v(t)<-pv(b)$. In particular, $pv(a)+v(t)$ lies in the convex hull of $v(b)$; but $v(t)-(pv(a)+v(t))$ is $p$-divisible, contradicting the fact that $v(b)\in F_p(v(t))$.

Finally, assume that $v(b)<0$ and $v(b)\notin F_p(v(t))$. This means that there is $h\in (pv(b),-pv(b))$ such that $v(t)+h$ is $p$-divisible. If $h>0$, we note that $v(t)+h-ph$ is also $p$-divisible, but $h-ph$ is now negative, so we might assume $h<0$. Now $v(t)+h<v(t)$ and is $p$-divisible. Take $a\in K$ such that $pv(a)=-(v(t)+h)$. Now $pv(a)>-v(t)$ so $a\in I_t$. But $p(v(ab))=pv(a)+pv(b)=-v(t)-h+pv(b)<-v(t)$ since $pv(b)<h$, thus $ab\notin I_t$, which in turn means $b\notin S_t$.
  \end{proof}

Now \Cref{str-abc} follows easily: in one direction, if $G\substruct G\oplus H$ with $H\neq\{0\}$, then we take $k=\mathbb C((\bigoplus_{\omega} H))$, and by \Cref{thm:straugvf}, since $k$ and $G$ are strongly co-augmentable, the valuation is not definable in $k((G))$.

In the other direction, if $G$ is not augmentable by infinitesimals, then by \Cref{thm:StronglyAugmentableInfinitesimalsOAGs}, there is $n\geqslant 2$ and $g\in G$ such that $F_n(g)=\{0\}$. By \Cref{lem:FnSChmitt}, there is $p\vert n$ and $h\in G$ such that $F_p(h)=\{0\}$. We conclude by \Cref{grf}.

\begin{remark}\label{rem:not-axomatisable2}
As in \Cref{rem:not-axomatisable}, we note that the class of ordered abelian groups which are strongly augmentable by infinitesimals is axiomatisable in $\Loag$, in fact, an axiomatization is explicitely given in \Cref{thm:StronglyAugmentableInfinitesimalsOAGs}.

Its complement, however, is not, as any finite collection of conditions of the form $(\forall g) F_n(g)\neq \{0\}$ is satisfied by $\mathbb Z_{(p)}$ for $p$ big enough -- but $\mathbb Z_{(p)}$ satisfies $F_p(1)=\{0\}$.

Therefore, the class of ordered abelian groups with automatic definability is not axiomatisable.
\end{remark}

Furthermore, let us note that in this case, we have an explicit formula defining the valuation ring, with quantifier alternation of the type $\forall\exists$.
Finally, let's compare our result to the main theorem of \cite{KKL}:

\begin{fact}[{\cite[Theorem 3.1]{KKL}}]
 Let $G$ be an ordered abelian group. If $G$ is not a closed subset of its divisible hull $\Gdiv$, then $G$ has automatic definability.
\end{fact}

\Cref{str-abc} gives a new proof of their result.
Moreover, we now show that the results of Kuhlmann, Krapp and Link
are essentially optimal:
%, we can give another proof of their result. In fact, we even prove the following:

\begin{theorem}
Let $G$ be an ordered abelian group. Then $G$ is strongly augmentable by infinitesimals if and only if $G$ has no minimum positive element and is closed in its divisible hull $\Gdiv$. \label{Gdiv}
\end{theorem} 

\begin{proof}
 We use the characterization of strong augmentability by infinitesimals given in \Cref{thm:StronglyAugmentableInfinitesimalsOAGs} as well as \cite[Remark 3.2]{KKL}, which states the class of  ordered abelian groups $G$ which are closed in their divisible hulls is axiomatised by the axiom scheme $$(\varphi_n)_{n\geqslant 2}\colon (\forall a) \left [(\forall b>0)(\exists g)~ (|a-ng|<b\rightarrow (\exists h) (a=nh) )\right].$$

 Assume first that $G$ is strongly augmentable by infinitesimals. It is clear that $G$ cannot have a minimum positive element. Now fix $n\geqslant 2$ and $a\in G$. We know that $F_n(a)\neq\{0\}$, and there is two possibilities: either $F_n(a)=\emptyset$, in which case $a$ is $n$-divisible; or $F_n(a)\supsetneq\{0\}$, in which case $a-ng\notin F_n(a)$ for any $g\in G$ by definition of $F_n$, thus, taking any $b\in F_n(a)\setminus\{0\}$, we have that $|a-ng|>b$ for any $g\in G$. In both cases, we have that $\varphi_n$ holds in $G$.

 Assume now that $G$ has no minimum positive element and is closed in $\Gdiv$. Fix $n\geqslant 2$ and  $a\in G$. If $a$ is $n$-divisible, then $F_n(a)=\emptyset$. If $a$ is not $n$-divisible, then by $\varphi_n$, there must exist some $b\in G_{>0}$ such that for all $g\in G$, $|a-ng|\geqslant b$. Since $G$ has no minimum positive element, we move if needed to some elementary extension $G^*\superstruct G$ with $A(b)\supsetneq\{0\}$ -- where $A(b)$ is the largest convex subgroup not containing $b$. We have by elementarity that for all $g\in G^*$, $|a-ng|\geqslant b$. In particular, $a-ng\notin A(b)$, and $F_n(a,G^*)\supseteq A(b)$, and thus $F_n(a,G)\neq\{0\}$ by elementarity. This holds for all $n \geqslant 2$, and therefore by \Cref{thm:StronglyAugmentableInfinitesimalsOAGs}, $G$ is strongly augmentable. 
\end{proof}

\subsubsection{Without parameters}

Similarly as for the strong case, if an ordered abelian group $G$ is weakly augmentable by infinitesimals, that is, if $G\equiv G\oplus H$ for some ordered abelian group $H\neq\{0\}$, then we can construct a henselian valued field $(K,v)$ with $vK\equiv G$ and on which $v$ is not $\emptyset$-definable. We now prove that this is the only obstacle to definability:

\begin{theorem}\label{weakABC}
 Let $G$ be an ordered abelian group. The following are equivalent:
 \begin{enumerate}
\item $G$ has automatic $\emptyset$-definability.
\item The valuation $v$ is uniformly definable in the class of henselian valued fields with value group elementary equivalent to $G$.
\item $G$ is not weakly augmentable by infinitesimals.
 \end{enumerate}
\end{theorem}

\begin{proof}
We already argued why if $G$ is weakly augmentable by infinitesimals, then it can't have automatic $\empty$-definability, that is, (1) implies (3). (2) implies (1) is clear, and to prove (3) implies (2), we apply Beth's definability theorem (\Cref{fact:beth}) with $\Lcal=\Lring$, $\Lcal'=\Lval$, $T$ the $\Lcal'$-theory of all henselian valued fields with value group elementary equivalent to $G$, and the $\mathcal{L}_\textrm{val}$-formula $\phi(x)$ which asserts that $ x \in \mathcal{O}_v$.

Let $(K,v_1)$ and $(K,v_2)$ be two models of $T$ which are expansions of the same field $K$, we want to prove that $\mathcal O_{v_1}=\mathcal O_{v_2}$.

\emph{Case 1: the residue fields $Kv_1$ and $Kv_2$ are both separably closed.} Then we have $\mathcal O_{v_1}\subseteq \mathcal O_{v_K}$, and $Kv_K$ itself is separably closed, by the defining property of the canonical henselian valuation. But if $\mathcal O_{v_1}$ is a strict refinement of $\mathcal O_{v_K}$, then $v_1$ induces a non-trivial valuation $\overline{v_1}$ on $Kv_K$. Since $Kv_K$ is separably closed, the value group $\overline {v_1}(Kv_K)$ is divisible. Now, any divisible ordered abelian group is weakly (even strongly) augmentable by infinitesimals, and therefore, $G\equiv v_1K\equiv v_KK\oplus\overline {v_1}(Kv_K)$ is also weakly-augmentable by infinitesimals, a contradiction. Thus $\mathcal O_{v_1}=\mathcal O_{v_K}$, and similarly, $\mathcal O_{v_2}=\mathcal O_{v_K}$, so they are indeed equal.

\emph{Case 2: at least one of the residue fields $Kv_1$ or $Kv_2$ is not separably closed.} Then at least one of $v_1$ or $v_2$ is a coarsening (not necessarily strict) of $v_K$, and is thusfore comparable with all other henselian valuations. So $v_1$ and $v_2$ are comparable, say $\mathcal O_{v_1}\subseteq\mathcal O_{v_2}$, and denote by $\overline{v_1}$ the induced valuation on $Kv_2$. But now $G\equiv v_1K\equiv v_2K\oplus\overline{v_1}(Kv_2)\equiv G\oplus\overline{v_1}(Kv_2)$. Because $G$ is not augmentable, we conclude that $\overline{v_1}(Kv_2)=\{0\}$, that is, $\overline{v_1}$ is trivial, or $v_1=v_2$.
\end{proof}

\begin{remark}\label{rem:not-axomatisable3}
As in \Cref{rem:not-axomatisable,rem:not-axomatisable2}, we have that the class of ordered abelian groups augmentable by infinitesimals is axiomatisable, but its complement is not.

Indeed, this class is clearly closed under elementary equivalence and ultraproducts, since $\prod_{\mathcal U}(G_i\oplus H_i)\cong\prod_{\mathcal U}(G_i)\oplus\prod_{\mathcal U}(H_i)$.

To see that its complement is not closed under ultraproducts, we simply observe that $\mathbb Z_{(p)}$ is not weakly augmentable by infinitesimals, but $\prod_{\mathcal U}\mathbb Z_{(p)}\equiv \mathbb Q$ is.
\end{remark}

In \Cref{prop:theory-of-weak-aug-lo}, we give an explicit axiomatisation of the theory of weakly augmentable coloured linear orders. From this, we can write the theory of ordered abelian groups with weakly augmentable $n$-spine for any fixed collection of $n\in\mathbb N$. This is not however enough to obtain an explicit axiomatisation of the class of ordered abelian groups weakly augmentable by infinitesimals.

Finally, we give an example on how \Cref{weakABC} can be used to give a new proof of previously know results. Note that by the sheer virtue of being an equivalence, \Cref{weakABC} is stronger than any theorem of the form ``if $(K,v)$ is a henselian valued field and $vK$ has property $P$, then $v$ is $\emptyset$-definable''.

\begin{fact}[\cite{Hong}]
 Let $G$ be an ordered abelian group. If $G$ is $p$-regular and not $p$-divisible for some prime $p$, then $G$ has automatic $\emptyset$-definability.
\end{fact}

\begin{proof}
 Let $G$ be such a group and suppose that $G\equiv G\oplus H$ for some $H\neq\{0\}$. Now $G\oplus H$ is $p$-regular by elementarity, and a quotient of a $p$-regular ordered abelian group by a non-trivial convex subgroup is always $p$-divisible, thus $G\equiv G\oplus H/H$ is $p$-divisible, a contradiction.
\end{proof}

For an example of a non-regular ordered abelian group which nonetheless has automatic $\emptyset$-definability, consider $\bigoplus_{p\in\mathbb P}\mathbb Z_{(p)}\oplus\mathbb Z$.
%auto-ignore
\section{Augmentability of Coloured Linear Orders}\label{sec:aug-lo}
In this section, we fix a language of coloured linear order $\lL \coloneqq \{<, \chi_\alpha : \alpha<\lambda\}$, where $\chi_\alpha$'s are all unary predicates. We call \emph{chain} any coloured linear order in this language. Recall that, given two chains $\mathcal X=(X,<_X)$ and $\mathcal Y=(Y,<_Y)$, we let $\mathcal X+\mathcal Y$ denote the chain with base set $X\sqcup Y$ and with ordering ``$<_+$'' defined as $x<_+ y$ if and only if $x,y\in X\wedge x<_X y$, or $x,y\in Y\wedge x<_Y y$, or $x\in X\wedge y\in Y$. We define colours in $\mathcal X+\mathcal Y$ in the obvious way.
%If $(X_y)_{y\in Y}$ are orders indexed by another order $Y$, we let $\sum_Y X_y$ denote the lexicographic sum, that is the set $\{(y,x) : y\in Y, x\in X_y\}$ with the lexicographic order. If all $X_y$ are copies of the same linear order $X$, we denote it $Y\cdot X$ or $\sum_Y X$.

\begin{definition}\label{def:augmentableLinearOrder}
Let $X$ and $Y$ be chains. A chain $X$ is said to be 
\begin{itemize}
    \item weakly augmentable (on the right) if there exists an non-empty chain $Y$ such that $X\equiv X+Y$. We say then that $X$ is weakly augmentable (on the right) by $Y$, and $Y$ is a \emph{weak augment of $X$}.
    \item strongly augmentable (on the right) if there exists an non-empty chain $Y$ such that $X\substruct X+Y$. We say then that $X$ is strongly augmentable (on the right) by $Y$, and $Y$ is a \emph{strong augment of $X$}.
\end{itemize}
\end{definition}

Similarly, we can define augmentability \emph{on the left}, by replacing $X+Y$ with $Y+X$. As the situation is symmetric, we mostly consider, in this section, augmentability on the right, and will only specify when needed.
It is clear that being  strongly augmentable implies being weakly augmentable.

\begin{example}\ 
\begin{itemize}
\item $\omega$ is strongly augmentable by $\zeta$.
\item $\omega^*$ is weakly augmentable by a single point. However, it is not strongly augmentable, as it has a maximum.
\item $\omega+1$ is not weakly augmentable.

\item $\sum_{\omega} (\zeta + n + \zeta)$ is strongly augmentable by $\sum_{\zeta} (\zeta + \omega +\omega^* + \zeta)$.

\item for any chain $X$, $\sum_{\omega^*}X$ is weakly augmentable by $X$.
\end{itemize}
\end{example}
The following lemma, although not used in this section, is useful in general to study augmentability. In particular, it was used in the proof of \Cref{thm:WeakAugmentableInfinitesimalsOags}.

\begin{lemma}\label{lem:StrongAugment}
    Let $X$ be a strongly augmentable chain, and let $Y$ be a strong augment of $X$:
    \[ X \substruct X+Y\]
    Then $\sum_\mathbb{\zeta} Y$ is also a strong augment of $X$:
    \[ X \substruct X+\sum_\mathbb{\zeta} Y.\]
\end{lemma}
\begin{proof}
    This is a direct application of Tarski's chain lemma \cite[Theorem~2.1.4 \& Exercise~2.1.1]{TZ12} :
    Define $L_n$ the chain consisting of $2n+1$ successive copies of $Y$:
    \[L_n = \sum_{-n}^nY  = Y_{-n} + \dots+ Y_n\]
    where $Y_i$ is a copy of $Y$ for all natural number $i$. With the obvious embedding of $L_n$ in $L_{n+1}$, we have for all $n\in \mathbb{N}$:
    \[X+L_n \substruct  X + Y_{-n-1} + L_n + Y_{n+1} = X+ L_{n+1}.\]
    By Tarski's chain lemma,  $X \substruct \bigcup_n (X+L_n) = X+\sum_\mathbb{\zeta} Y$, as wanted.
    \end{proof}

Let $A \subseteq B$ be two chains, and $a\in B$ an element. We denote:
\begin{itemize}
    \item $B_{\square A} = \{b\in B : \exists a\in A,~  b\square a \}$ for $\square\in \{\geqslant,\leqslant\}$,
    \item $B_{\square A} = \{b\in B : \forall a\in A,~  b \square a \}$ for $\square\in \{>,<\}$,
        \item $\langle A\rangle_{B} \coloneqq B_{\geqslant A} \cap B_{\leqslant A} = \{b\in B : \exists a,a'\in A,~   a \leqslant b \leqslant a' \}$,
    \item $B_{\square a} \coloneqq B_{\square \{a\}}$ for $\square\in \{\geqslant, >, \leqslant, <\}$.
    
\end{itemize}

\begin{lemma}[See {\cite[Example 1]{DelonLucas}}]\label{lem:hulls!}
 Let $A,B$ be chains with $A\substruct B$. Then:
 \begin{enumerate}
  \item $A\substruct \langle A\rangle_B\substruct B$,
  \item $A\substruct B_{\leqslant A}\substruct B$, and
  \item $A\substruct B_{\geqslant A}\substruct B$.
 \end{enumerate}
\end{lemma}

\begin{proof}
 Let $a\in A$. With the notation above, we have $B_{\leqslant a}=\sset{x\in B}{x\leqslant a}$ and $A_{>a}=\sset{x\in A}{x>a}$. Denote $B_{\leqslant a}+A_{>a}$ by $C_a$. We have naturally $A_{\leqslant a}\substruct B_{\leqslant a}$, therefore, $A=A_{\leqslant a}+A_{>a}\substruct C_a\substruct B$. The chain of structures $(C_a)_{a\in A}$ is thus elementary, and by Tarski's chain lemma, $A\substruct \bigcup_{a\in A} C_a\substruct B$, and $\bigcup_{a\in A} C_a$ is exactly $B_{\leqslant A}$.
 
 Reversing the ordering gives a proof of $A\substruct B_{\geqslant A}\substruct B$. For the convex hull, we simply note that ${(B_{\leqslant A})}_{\geqslant A}=\langle A\rangle_B$.
\end{proof}

We now use this fact to characterize strongly augmentable linear orders:

\begin{corollary}\label{lem:aug-mo}
 Let $A$ be a chain. Then $A$ is strongly augmentable on the right (resp.~on the left) if and only if it has no last element (resp.~no first element).
\end{corollary}

\begin{proof}
We do the proof on the right; reversing the ordering gives a proof on the left. If $A$ has a last element, clearly there can't be a non-empty $X$ such that $A\substruct A+X$. Suppose that $A$ has no last element. The ``type at infinity'' 
\[\pi_{\infty}(x)=\sset{a\leqslant x}{a\in A}\]
is thus satisfiable.
Let $B\superstruct A$ be $|A|^+$-saturated. Then in particular, $B$ contains a realisation of $\pi_{\infty}$, that is, $B_{>A}\neq\emptyset$.

\begin{remark}
 For the first (and last) time in this paper, we have a class which is finitely axiomatisable, as ``not having a max'' is expressed by one sentence. This means the classes of strongly augmentable and not strongly augmentable chains are both axiomatisable.
\end{remark}

Now by \Cref{lem:hulls!}, we have $A\substruct B_{\leqslant A}$, therefore by back-and-forth, we have $A+X\substruct B_{\leqslant A}+X$ for any coloured linear order $X$; in particular, $A+B_{>A}\substruct B_{\leqslant A}+B_{>A}=B$. Now $A\substruct B$ and $A\subset A+B_{>A}\substruct B$, therefore, $A\substruct A+B_{>A}$.
\end{proof}

The following proposition, also due to Delon and Lucas  (see \cite[p.510]{DelonLucas}), gives a simple characterisation of weakly augmentable chains.

\begin{proposition}\label{prop:ChainAugmentable} 
    Let $L$ be a non-empty chain. A chain $C$ is weakly-augmentable by $L$ if and only if there is a chain $B$ such that $C \equiv B+ \sum_{\omega^*} L$.
\end{proposition}
\begin{proof}
    Clearly, if $C$ is of the form $B+ \sum_{\omega^*} L$, it is weakly augmentable by $L$. Assume now that $C$ is weakly augmentable by $L$. We may replace $C$ with a saturated enough model, so that there is an isomorphism $\sigma: C+L \rightarrow C$. The application $\sigma$ can then be iterated on $L$, since $\sigma(L)$ is a subset of $C$. Letting $L^n$ denote $\sigma^{(n)}(L)$, all $L^n$ are clearly disjoint, with $L^{n+1}<L^{n}$. Their union $ \bigcup_n L^n = \sum_{\omega^*} L^n$ is an end-segment of $C$. Denote $B:= C \setminus \sum_{\omega^*} L^n$. Then by back-and-forth, we have
     \[C=B  +  \sum_{\omega^*} L^n \equiv B + \sum_{\omega^*} L ,\]
     and this conclude the proof.
\end{proof}

This characterisation of weakly augmentable chain is similar to the characterisation of ``right-abosrbing'' chains in \cite[Proposition 3.15]{EP}: a chain $X$ is called right-absorbing if there exists $Y$ such that $X+Y\simeq X$, a property that we could instead call ``isomorphically augmentable on the right''. \cite[Proposition 3.15]{EP} states that $X$ is right-absorbing if and only if $X\simeq X+\sum_{\omega^*}A$ for some chain $A$.

To conclude this short section, we note that the class of weakly augmentable chains is also elementary, and an explicit axiomatisation can be given. We are deeply grateful to David Gonzalez for providing us such an axiomatisation, which we now present. 
Let $\varphi$ be a sentence in $\lL$, and let $x$ be a new variable which doesn't occur in $\varphi$. We denote by $\varphi_{\leqslant}(x)$ the formula with free variable $x$ obtained from $\varphi$ by replacing quantifiers $(\exists y)$ and $(\forall y)$ by $(\exists y\leqslant x)$ and $(\forall y\leqslant x)$ respectively.
    It is easy to see that for all chains $C$ and $c \in C$, we have:
    \[ C \models \varphi_\leqslant(c) \Leftrightarrow C_{\leqslant c} \models \varphi. \]

 As a preliminary remark, note that if a chain does not have a maximal element, then it is strongly augmentable by \Cref{lem:aug-mo} -- and thus also weakly augmentable. Remains to axiomatise weakly augmentable chains with a maximal element. For ease of notation, we move to a language containing a new constant symbol standing for the maximal element.
    
\begin{proposition}\label{prop:theory-of-weak-aug-lo}
    Let $\lL_M \coloneqq \lL \cup \{M\}$ be the language $\lL$ with a constant $M$. 
    The theory $T_w$ of weakly augmentable chains with maximal element $M$ consists of the following axioms:
        \begin{itemize}
            \item $ (\forall y) ~(y\leqslant M)$, and
            \item for each sentence $\varphi$ in $\lL$, ($\varphi \rightarrow (\exists x<M) ~ \varphi_\leqslant (x)$).
        \end{itemize}
\end{proposition}

\begin{proof}
    Let $C$ be a weakly augmentable chain with maximal element $M$. Then $C\equiv C+L$ for some non-empty $L$, and by elementarity, $L$ has a maximal element. Let $\vphi$ be a sentence in $\mathcal L$ such that $C\vDash\vphi$. Now $C=(C+L)_{\leqslant M_C}$, where $M_C$ is the maximal element of $C$. Thus, $C+L\vDash(\exists x<M)\vphi_\leqslant(x)$, and $C$ alike.
    
    Conversely, assume that $C$ is a model of $T_w$. Without loss of generality, we may assume that $C$ is $|\lL|^+$-saturated. By compactness, we can find an element $m\in C$ such that $m<M$ and $C_{\leqslant m} \equiv C$. Then $C = C_{\leqslant m}+C_{> m} \equiv C+C_{> m}$. In particular, $C$ is weakly augmentable. 
\end{proof}

\begin{remark}
 The class of weakly augmentable chains is not finitely axiomatisable, since its complement is not axiomatisable: any finite chain $n$ is not weakly augmentable, but for a non-principal ultrafilter $\mathcal U$, we have that $\prod_{\mathcal U}n\simeq\omega+\omega^*$ is weakly augmentable. 
\end{remark}

%\listoftodos{}

\bibliographystyle{abbrv}
\bibliography{bibtex}

\end{document}